\documentclass[a4paper,english]{article}
\usepackage[T1]{fontenc}
\usepackage[latin9]{inputenc}
\usepackage{color}
\usepackage{url}
\usepackage{amsthm}
\usepackage{amsmath}
\usepackage{graphicx}
\usepackage{amssymb}

\theoremstyle{plain}

\usepackage{a4wide}
\usepackage{color}
\usepackage{graphics}
\usepackage{floatflt}
\usepackage{amsmath, amsthm, amssymb}

\usepackage{babel}

\begin{document}
\newcommand{\np}{\newpage}

\newcommand{\vs}{\vspace{1cm}}

\newtheorem{prop}{Proposition} 
\newtheorem{theo}{Theorem} 
\newtheorem{coro}{Corollary} 
\newtheorem{lem}{Lemma} 
\newtheorem{conj}{Conjecture}
\newtheorem{Def}{Definition}

\newcommand{\bp}{\begin{prop}}
\newcommand{\ep}{\end{prop}}

\title{Semiclassical origin of the spectral gap for transfer operators of
partially expanding map.}

\author{Frédéric Faure\textit{}%
\thanks{Institut Fourier, 100 rue des Maths, BP74 38402 St Martin d'Hères.
frederic.faure@ujf-grenoble.fr \protect\url{http://www-fourier.ujf-grenoble.fr/~faure}%
}}
\maketitle
\begin{abstract}
We consider a simple model of partially expanding map on the torus.
We study the spectrum of the Ruelle transfer operator and show that
in the limit of high frequencies in the neutral direction (this is
a semiclassical limit), the spectrum develops a spectral gap, for
a generic map. This result has already been obtained by M. Tsujii
in \cite{tsujii_05}. The novelty here is that we use semiclassical
analysis which provides a different and quite natural description.
We show that the transfer operator is a semiclassical operator with
a well defined {}``classical dynamics'' on the cotangent space.
This classical dynamics has a {}``trapped set'' which is responsible
for the Ruelle resonances spectrum. In particular we show that the
spectral gap is closely related to a specific dynamical property of
this trapped set.
\end{abstract}
\footnote{2000 Mathematics Subject Classification:37D20 hyperbolic systems (expanding,
Anosov, Axiom A, etc.) 37C30 Zeta functions, (Ruelle-Frobenius) transfer
operators, and other functional analytic techniques in dynamical systems
81Q20 Semi-classical techniques

Keywords: Transfer operator, Ruelle resonances, decay of correlations,
Semi-classical analysis. %
}

\tableofcontents{}

\section{Introduction}

Chaotic behavior of certain dynamical systems is due to hyperbolicity
of the trajectories. This means that the trajectories of two closed
initial points will diverge from each other either in the future or
in the past (or both). As a result the behavior of an individual trajectory
looks like complicated and unpredictable. However evolution of a cloud
of points seems more simple: it will spread and equidistributes according
to an invariant measure, called an equilibrium measure (or S.R.B.
measure). Following this idea, D. Ruelle in the 70' \cite{ruelle_75,ruelle_86},
has shown that instead of considering individual trajectories, it
is much more natural to consider evolution of densities under a linear
operator called the Ruelle Transfer operator or the Perron Frobenius
operator.

For dynamical systems with strong chaotic properties, such as uniformly
expanding maps or uniformly hyperbolic maps, Ruelle, Bowen, Fried,
Rugh and others, using symbolic dynamics techniques, have shown that
the spectrum of the transfer operator has a discrete spectrum of eigenvalues.
This spectral description has an important meaning for the dynamics
since each eigenvector corresponds to an invariant distribution (up
to a time factor). From this spectral characterization of the transfer
operator, one can derive other specific properties of the dynamics
such as decay of time correlation functions, central limit theorem,
mixing, ... In particular a spectral gap implies exponential decay
of correlations. 

This spectral approach has recently (2002-2005) been improved by M.
Blank, S. Gouëzel, G. Keller, C. Liverani \cite{liverani_02,liverani_04,liverani_05}
and V. Baladi and M. Tsujii \cite{baladi_sobolev_05,baladi_05} (see
\cite{baladi_05} for some historical remarks), through the construction
of functional spaces adapted to the dynamics, independent of any symbolic
dynamics. The case of dynamical systems with continuous time is more
delicate (see \cite{melbourne_07} for historical remarks). This is
due to the direction of time flow which is neutral (i.e. two nearby
points on the same trajectory will not diverge from each other). In
1998 Dolgopyat \cite{dolgopyat_98,dolgopyat_02} showed the exponential
decay of correlation functions for certain Anosov flows, using techniques
of oscillatory integrals and symbolic dynamics. In 2004 Liverani \cite{liverani_contact_04}
adapted Dolgopyat's ideas to his functional analytic approach, to
treat the case of contact Anosov flows. In 2005 M. Tsujii \cite{tsujii_05}
obtained an explicit estimate for the spectral gap for the suspension
of an expanding map. Then in 2008 M. Tsujii \cite{tsujii_08} obtained
an explicit estimate for the spectral gap, in the case of contact
Anosov flows.

\paragraph{Semiclassical approach for transfer operators:}

It also appeared recently \cite{fred-RP-06,fred-roy-sjostrand-07}
that for hyperbolic dynamics, the study of transfer operator is naturally
a semiclassical problem in the sense that a transfer operator can
be considered as a {}``Fourier integral operator'' and using standard
tools of semiclassical analysis, some of its spectral properties can
be obtained from the study of {}``the associated classical symplectic
dynamics'', namely the initial hyperbolic dynamics lifted on the
cotangent space (the phase space).

The simple idea behind this, crudely speaking, is that a transfer
operator transports a {}``wave packet'' (i.e. localized both in
space and in Fourier space) into another wave packet, and this is
exactly the characterization of a Fourier integral operator. A wave
packet is characterized by a point in phase space (its position and
its momentum), hence one is naturally led to study the dynamics in
phase space. Moreover, since any function or distribution can be decomposed
as a linear superposition of wave packets, the dynamics of wave packets
characterizes completely the transfer operators.

Following this approach, in the papers \cite{fred-RP-06,fred-roy-sjostrand-07}
we studied hyperbolic diffeomorphisms. The aim of the present paper
is to show that semiclassical analysis is also well adapted (even
better) for hyperbolic systems with neutral direction. We consider
here the simplest model: a partially expanding map $f:\left(x,s\right)\rightarrow f\left(x,s\right)$,
i.e. a map on a torus $\left(x,s\right)\in S^{1}\times S^{1}$ with
an expanding direction $\left(x\in S_{x}^{1}\right)$ and a neutral
direction $\left(s\in S_{s}^{1}\right)$ (the inverse map $f^{-1}$
is $k$-valued, with $k\geq2$). The results are presented in section
2. We summarize them in few lines. First in order to reduce the problem
and drop out the neutral direction, we use a Fourier analysis in $s\in S_{s}^{1}$
and decompose the transfer operator $\hat{F}$ on $S_{x}^{1}\times S_{s}^{1}$
(defined by $\hat{F}\varphi:=\varphi\circ f$) as a collection of
transfer operators $\hat{F}_{\nu}$ on the expanding space $S_{x}^{1}$
only, with $\nu\in\mathbb{Z}$ being the Fourier parameter and playing
the role of the semiclassical parameter. The semiclassical limit is
$\left|\nu\right|\rightarrow\infty$.

Then we introduce a (multivalued) map $F_{\nu}$ on the cotangent
space $\left(x,\xi\right)\in T^{*}S_{x}^{1}$ which is the canonical
map associated to the transfer operator $\hat{F}_{\nu}$. The fact
that the initial map $f$ is expanding along the space $S_{x}^{1}$
implies that on the cylinder $T^{*}S_{x}^{1}$ trajectories starting
from a large enough value of $\left|\xi\right|$ escape towards infinity
($\left|\xi\right|\rightarrow\infty$). We define the trapped set
as the compact set $K=\lim_{n\rightarrow\infty}F_{\nu}^{-n}\left(K_{0}\right)$
where $K_{0}\subset T^{*}S^{1}$ is an initial large compact set.
$K$ contains trajectories which do not escape towards infinity. 

Using a standard semiclassical approach (with escape functions on
phase space \cite{sjostrand_87}) we first show that the operator
$\hat{F}_{\nu}$ as a discrete spectrum called Ruelle resonances (we
have to consider $\hat{F}_{\nu}$ in Sobolev space of distributions).
This is Theorem \ref{th:spectrum_resonances}. This result is well
known, but the semiclassical approach we use here is new.

Then we show that a specific hypothesis on the trapped set implies
that the operator $\hat{F}_{\nu}$ develops a {}``spectral gap''
in the semi-classical limit $\nu\rightarrow\infty$ (i.e. its spectral
radius reduces). This is Theorem \ref{th:gap_spectral} illustrated
on Figure \ref{fig:Spectrum-of-F}. This Theorem is very similar to
Theorem 1.1 in \cite{tsujii_05}. With the semiclassical approach,
this result is very intuitive: the basic idea (followed in the proof)
is that an initial wave packet $\varphi_{0}$ represented as a point
on the trapped set $K$ evolves in several wave packets $\left(\varphi_{j}\right)_{j=1\rightarrow k}$under
the transfer operator $\hat{F}_{\nu}$, but in general only one wave
packet remains on the trapped set $K$ and the $\left(k-1\right)$
other ones escape towards infinity. As a result the probability on
the trapped set $K$ decays by a factor $1/k$. This is the origin
of the spectral gap at $1/\sqrt{k}$ on Figure \ref{fig:Spectrum-of-F}.

This work has been supported by {}``Agence Nationale de la Recherche''
under the grant JC05\_52556.

\section{Model and results}

\subsection{A partially expanding map}

Let $g:S^{1}\rightarrow S^{1}$ be a $C^{\infty}$ diffeomorphism
(on $S^{1}:=\mathbb{R}/\mathbb{Z}$). $g$ can be written as $g:\mathbb{R}\rightarrow\mathbb{R}$
with $g\left(x+1\right)=g\left(x\right)+1,\,\forall x\in\mathbb{R}$.
Let $k\in\mathbb{N}$, $k\geq2$, and let the map $E:S^{1}\rightarrow S^{1}$
be defined by\begin{equation}
E:\quad x\in S^{1}\rightarrow E\left(x\right)=kg\left(x\right)\, mod\,1\label{eq:def_E}\end{equation}
Let \[
E_{min}:=\min_{x}\left(\frac{dE}{dx}\right)\left(x\right)=k\min_{x}\left(\frac{dg}{dx}\left(x\right)\right)\]
We will suppose that the function $g$ is such that \begin{equation}
\boxed{E_{min}>1}\label{eq:Emin}\end{equation}
 so that $E$ is a\textbf{ uniform expanding map} on $S^{1}$. The
map $E$ is then a $k:1$ map (i.e. every point $y$ has $k$ previous
images $x\in E^{-1}\left(y\right)$). Let $\tau:S^{1}\rightarrow\mathbb{R}$
be a $C^{\infty}$ function, and define a map $f$ on $\mathbb{T}^{2}=S^{1}\times S^{1}$
by:

\vspace{0.0cm}\begin{center}{\color{red}\fbox{\color{black}\parbox{15cm}{\begin{equation}
f:\quad\left(\begin{array}{c}
x\\
s\end{array}\right)\longmapsto\left(\begin{array}{l}
x'=E\left(x\right)=kg\left(x\right)\qquad mod\,1\\
s'=s+\frac{1}{2\pi}\tau\left(x\right)\qquad mod\,1\end{array}\right)\label{eq:def_f}\end{equation}

}}}\end{center}\vspace{0.0cm}

\begin{figure}[h]
\begin{centering}
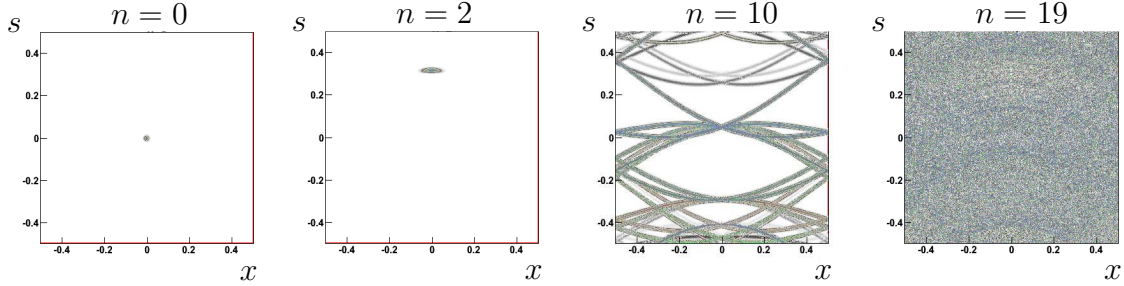
\par\end{centering}

\centering{}\caption{\label{fig:Evolution}Numerical evolution of an initial small cloud
of points on the torus $\left(x,s\right)\in\mathbb{T}^{2}$ under
the map $f$, Eq.(\ref{eq:def_f}), at different time $n=0,2,10,19$.
We have chosen here $E\left(x\right)=2x$ and $\tau\left(x\right)=\cos\left(2\pi x\right)$.
The initial cloud of points is centered around the point $\left(0,0\right)$.
For small time $n$, the cloud of point is transported in the vertical
direction $s$ and spreads in the expanding horizontal direction $x$.
Due to instability in $x$ and periodicity, the cloud fills the torus
$S^{1}\times S^{1}$ for large time $n$. On the last image $n=19$,
one observes an invariant absolutely continuous probability measure
(called SRB measure, equal to the Lebesgue measure in our example).
It reveals the mixing property of the map $f$ in this example.}

\end{figure}

The map $f$ is also a $k:1$ map. The map $f$ is a very simple example
of a \textbf{compact group extension} of the expanding map $E$ (see
\cite{dolgopyat_02}, \cite[p.17]{pesin_04}). It is also a special
example of a partially hyperbolic map%
\footnote{A even more general setting would be a $C^{\infty}$ map $f:M\rightarrow M$
on a compact Riemannian manifold $M$, which is supposed to be partially
expanding, i.e., for any $m\in M$, the tangent space $T_{m}M$ decomposes
continuously as\[
T_{m}M=E_{u}\left(m\right)\oplus E_{0}\left(m\right)\]
where $E_{u}\left(m\right)$ is a (non invariant) expanding direction
(with respect to a Riemannian metric $g$): \[
\left|D_{m}f\left(v_{u}\right)\right|_{g}>\left|v_{u}\right|_{g},\qquad\forall v_{u}\in E_{u}\left(m\right)\]
and $E_{0}\left(m\right)$ the neutral direction: there exist a non
zero global section $v_{0}\in C^{\infty}\left(TM\right)$ such that
$v_{0}\left(m\right)\in E_{0}\left(m\right)$ and $Df\left(v_{0}\right)=v_{0}$.
In our example (\ref{eq:def_f}), $M=S^{1}\times S^{1}$, the neutral
section is  $v_{0}=\left(0,1\right)$, and the expanding direction
$E_{u}\left(m\right)$ is spanned by the vector $\left(1,0\right)$.%
}. See figure \ref{fig:Evolution}.

\subsection{Transfer operator}

Instead of studying individual trajectories which have chaotic behavior,
one prefer to study the evolutions of densities induced by the map
$f$. This is the role of the \textbf{Perron-Frobenius transfer operator
$\hat{F}^{*}$ }on $C^{\infty}\left(\mathbb{T}^{2}\right)$ given
by:\begin{equation}
\left(\hat{F}^{*}\psi\right)\left(\mathbf{y}\right)=\sum_{\mathbf{x}\in f^{-1}\left(\mathbf{y}\right)}\frac{1}{\left|D_{\mathbf{x}}f\right|}\psi\left(\mathbf{x}\right),\qquad\psi\in C^{\infty}\left(\mathbb{T}^{2}\right).\label{eq:def_F*}\end{equation}

Indeed if the the function $\psi$ has its support in a vicinity of
$x$ then the support of $\hat{F}^{*}\psi$ is in a vicinity of $y=f\left(x\right)$.
To explain the Jacobian in the prefactor, one checks%
\footnote{Since $\mathbf{y}=f\left(\boldsymbol{x}\right)$, then $d\mathbf{y}=\left|D_{\mathbf{x}}f\right|d\mathbf{x}$,
and \[
\int_{\mathbb{T}^{2}}\left(\hat{F}^{*}\psi\right)\left(\mathbf{y}\right)d\mathbf{y}=\sum_{\mathbf{x}\in f^{-1}\left(\mathbf{y}\right)}\int_{\mathbb{T}^{2}}\frac{1}{\left|D_{\mathbf{x}}f\right|}\psi\left(\mathbf{x}\right)\left|D_{\mathbf{x}}f\right|d\mathbf{x}=\int_{\mathbb{T}^{2}}\psi\left(\mathbf{x}\right)d\mathbf{x}\]
} that $\int_{\mathbb{T}^{2}}\left(\hat{F}^{*}\psi\right)\left(\mathbf{y}\right)d\mathbf{y}=\int_{\mathbb{T}^{2}}\psi\left(\mathbf{x}\right)d\mathbf{x}$,
i.e. the total measure is preserved. 

The operator $\hat{F}^{*}$ extends to a bounded operator on $L^{2}\left(\mathbb{T}^{2},d\mathbf{x}\right)$.
Its $L^{2}$-adjoint written $\hat{F}$ is defined by $\left(\hat{F}^{*}\psi,\varphi\right)_{L^{2}}=\left(\psi,\hat{F}\varphi\right)_{L^{2}}$,
with the scalar product $\left(\psi,\varphi\right)_{L^{2}}:=\int_{\mathbb{T}^{2}}\overline{\psi}\left(\mathbf{x}\right)\varphi\left(\mathbf{x}\right)d\mathbf{x}$.
One checks easily that $\hat{F}$ has a simpler expression than $\hat{F}^{*}$:
it is the \textbf{pull back operator}, also called the \textbf{Koopman
operator, }or \textbf{Ruelle transfer operator} and given by: \begin{equation}
\boxed{\left(\hat{F}\psi\right)\left(\boldsymbol{x}\right)=\psi\left(f\left(\boldsymbol{x}\right)\right)}\label{eq:def_F_pull_back}\end{equation}

\subsection{The reduced transfer operator}

The particular form of the map (\ref{eq:def_f}) allows some simplifications.
Observe that for a function of the form \[
\psi\left(x,s\right)=\varphi\left(x\right)e^{i2\pi\nu s}\]
with $\nu\in\mathbb{Z}$ (i.e. a Fourier mode in $s$), then\[
\left(\hat{F}\psi\right)\left(x,s\right)=\varphi\left(E\left(x\right)\right)e^{i\nu\tau\left(x\right)}e^{i2\pi\nu s}.\]

Therefore the operator $\hat{F}$ preserves the following decomposition
in Fourier modes: \begin{equation}
L^{2}\left(\mathbb{T}^{2}\right)=\bigoplus_{\nu\in\mathbb{Z}}\mathcal{H}_{\nu},\qquad\qquad\mathcal{H}_{\nu}:=\left\{ \varphi\left(x\right)e^{i2\pi\nu s},\qquad\varphi\in L^{2}\left(S^{1}\right)\right\} \label{eq:decomp_L2_T2}\end{equation}
The space $\mathcal{H}_{\nu}$ and $L^{2}\left(S^{1}\right)$ are
unitary equivalent. For $\nu\in\mathbb{Z}$ given, the operator $\hat{F}$
restricted to the space $\mathcal{H}_{\nu}\equiv L^{2}\left(S^{1}\right)$,
written $\hat{F}_{\nu}$ is%
\footnote{Notice that the operator $\hat{F}_{\nu}$ appears to be a transfer
operator for the expanding map $E$ with an additional weight function
$e^{i\nu\tau\left(x\right)}$.%
}:\begin{equation}
\boxed{\left(\hat{F}_{\nu}\varphi\right)\left(x\right):=\varphi\left(E\left(x\right)\right)e^{i\nu\tau\left(x\right)}\qquad\varphi\in L^{2}\left(S^{1}\right)\equiv\mathcal{H}_{\nu}}\label{eq:def_F_hat}\end{equation}
and with respect to the orthogonal decomposition (\ref{eq:decomp_L2_T2}),
we can write:\[
\boxed{\hat{F}=\bigoplus_{\nu\in\mathbb{Z}}\hat{F}_{\nu}}\]
We will study the spectrum of this family of operators $\hat{F}_{\nu}$,
with parameter $\nu\in\mathbb{Z}$, and consider more generally a
real parameter $\nu\in\mathbb{R}$. We will see that the parameter
$\nu$ is a \textbf{semiclassical parameter}, and $\nu\rightarrow\infty$
is the \textbf{semiclassical limit}. (if $\nu\neq0$, $\nu=1/\hbar$
in usual notations \cite{martinez-01}).

\paragraph{Remarks:}
\begin{itemize}
\item For $\nu=0$, $\hat{F}_{0}$ has an obvious eigenfunction $\varphi\left(x\right)=1$,
with eigenvalue $1$. Except in special cases (e.g. $\tau=0$), there
is no other obvious eigenvalues for $\hat{F}_{\nu}$ in $L^{2}\left(S^{1}\right)$.
\end{itemize}

\subsection{Main results on the spectrum of the transfer operator $\hat{F}_{\nu}$}

We first observe that by duality, the operator $\hat{F}_{\nu}$ defined
in (\ref{eq:def_F_hat}) extends to the distribution space $\mathcal{D}'\left(S^{1}\right)$:\[
\hat{F}_{\nu}\left(\alpha\right)\left(\varphi\right)=\alpha\left(\hat{F}_{\nu}^{*}\left(\varphi\right)\right),\qquad\alpha\in\mathcal{D}'\left(S^{1}\right),\varphi\in C^{\infty}\left(S^{1}\right),\]

where the $L^{2}$-adjoint  $\hat{F}_{\nu}^{*}$ is given by \begin{equation}
\left(\hat{F}_{\nu}^{*}\varphi\right)\left(y\right)=\sum_{x\in E^{-1}\left(y\right)}\frac{e^{-i\nu\tau\left(x\right)}}{E'\left(x\right)}\varphi\left(x\right),\qquad\varphi\in C^{\infty}\left(S^{1}\right).\label{eq:F_nu_*}\end{equation}

Before giving the main results, remind that for $m\in\mathbb{R}$,
the Sobolev space $H^{m}\left(S^{1}\right)\subset\mathcal{D}'\left(S^{1}\right)$
consists in distributions (or continuous functions if $m>1/2$) such
that their Fourier series $\hat{\psi}\left(\xi\right)$ satisfy $\sum_{\xi\in\mathbb{Z}}\left|\left\langle \xi\right\rangle ^{m}\hat{\psi}\left(\xi\right)\right|^{2}<\infty$,
with $\left\langle \xi\right\rangle :=\left(1+\xi^{2}\right)^{1/2}$.
It can equivalently be written (\cite{taylor_tome1} p.271).\[
H^{m}\left(S^{1}\right):=\left\langle \hat{\xi}\right\rangle ^{-m}\left(L^{2}\left(S^{1}\right)\right)\]
with the differential operator $\hat{\xi}:=-i\frac{d}{dx}$. 

The following theorem is well known \cite{ruelle_86}. We will however
provide a new proof based on semiclassical analysis.

\vspace{0.cm}\begin{center}{\color{blue}\fbox{\color{black}\parbox{15cm}{
\begin{theo}\label{th:spectrum_resonances}\textbf{Discrete spectrum of resonances.}

Let $m<0$. The operator $\hat{F}_{\nu}$ leaves the Sobolev space
$H^{m}\left(S^{1}\right)$ invariant, and \[
\hat{F}_{\nu}:H^{m}\left(S^{1}\right)\rightarrow H^{m}\left(S^{1}\right)\]
is a bounded operator and can be written\begin{equation}
\hat{F}_{\nu}=\hat{R}+\hat{K}\label{eq:decomp_F_nu}\end{equation}
where $\hat{K}$ is a compact operator, and $\hat{R}$ has a small
norm: \begin{equation}
\left\Vert \hat{R}\right\Vert \leq r_{m}:=\frac{1}{E_{min}^{\left|m\right|}}\sqrt{\frac{k}{E_{min}}}.\label{eq:bound_R}\end{equation}
(the interesting situation is $m\ll0$, since the norm $\left\Vert \hat{R}\right\Vert $
shrinks to zero for $m\rightarrow-\infty$). 

Therefore, $\hat{F}_{\nu}$ has an essential spectral radius less
than $r_{m}$, which means that $\hat{F}_{\nu}$ has discrete (eventually
empty) spectrum of generalized eigenvalues $\lambda_{i}$ outside
the circle of radius $r_{m}$ (see \cite[prop. 6.9 p.499]{taylor_tome1}).
The eigenvalues $\lambda_{i}$ are called \textbf{Ruelle resonances}.
Together with their associated eigenspace, they do not depend on $m$
and are intrinsic to the transfer operator $\hat{F}_{\nu}$.

\end{theo}
}}}\end{center}\vspace{0.cm}

The following theorem is analogous to Theorem 1.1 in \cite{tsujii_05}.
However the approach and the proof we propose are different and rely
on semiclassical analysis.

\vspace{0.cm}\begin{center}{\color{blue}\fbox{\color{black}\parbox{15cm}{
\begin{theo}\label{th:gap_spectral}\textbf{Spectral gap in the semiclassical
limit.}

if the map $f$ is \textbf{partially captive} (definition given page
\pageref{def:partially_captive}) (and $m$ small enough), then the
spectral radius of the operator $\hat{F}_{\nu}:H^{m}\left(S^{1}\right)\rightarrow H^{m}\left(S^{1}\right)$
does not depend on $m$ and satisfies in the semi-classical limit
$\nu\rightarrow\infty$:\begin{equation}
r_{s}\left(\hat{F}_{\nu}\right)\leq\frac{1}{\sqrt{E_{min}}}+o\left(1\right)\label{eq:upper_bound}\end{equation}

which is strictly smaller than $1$ from (\pageref{eq:Emin}).

\end{theo}
}}}\end{center}\vspace{0.cm}

\begin{figure}[h]
\begin{centering}
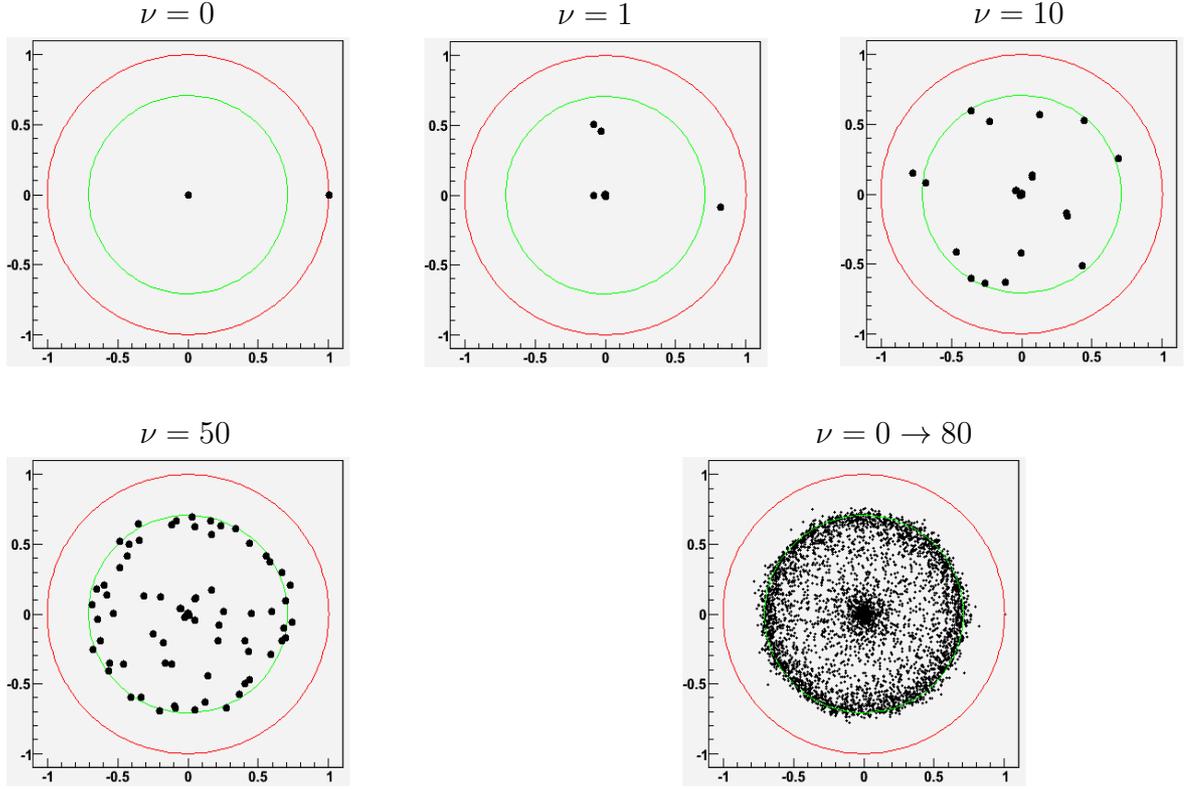
\par\end{centering}

\caption{\label{fig:Spectrum-of-F}Black dots are numerical computation of
the eigenvalues $\lambda_{i}$ of $\hat{F}_{\nu}$ for different values
of $\nu\in\mathbb{N}$, and union of these in the last image. We have
chosen here $E\left(x\right)=2x$ i.e. $k=2$, and $\tau\left(x\right)=\cos\left(2\pi x\right)$.
The external red circle has radius $1$. The internal green circle
has radius $1/\sqrt{E_{min}}=1/\sqrt{2}$ and represents the upper
bound given in Eq.(\ref{eq:upper_bound}). As $\nu\in\mathbb{R}$
moves continuously, the resonances move in a spectacular way. This
can be seen on a movie on \protect\url{http://www-fourier.ujf-grenoble.fr/~faure/articles}}

\end{figure}

\paragraph{Remarks:}
\begin{itemize}
\item This remark concerns the regularity of the eigenfunctions of $\hat{F}_{\nu}$.
Let $\lambda_{i}$ be a generalized eigenvalue of $\hat{F}_{\nu}$.
Let $\varphi_{i}$ denotes a generalized eigenfunction of $\hat{F}_{\nu}$
associated to $\lambda_{i}$ (i.e. $\hat{F}_{\nu}\varphi_{i}=\lambda_{i}\varphi_{i}$
if $\lambda_{i}$ is an eigenvalue). Then $\varphi_{i}$ belongs to
$H^{m}$ for any $m$ such that $m<m_{0}$ where $m_{0}$ is given
by $r_{m_{0}}=\left|\lambda_{i}\right|$. 
\item By duality we have similar spectral results for the Perron Frobenius
operator $\hat{F}_{\nu}^{*}:H^{m}\left(S^{1}\right)\rightarrow H^{m}\left(S^{1}\right)$
if $m>0$. The eigenvalues of $\hat{F}_{\nu}^{*}$ are $\overline{\lambda}_{i}$.
We have seen that the generalized eigenfunctions of $\hat{F}_{\nu}$
belong to different Sobolev spaces $H^{m}\left(S^{1}\right)$ with
$m<0$. Eq.(\ref{eq:bound_R}) says that $m$ should satisfy $r_{m}<\left|\lambda_{i}\right|$,
so $m\rightarrow-\infty$ as $\left|\lambda_{i}\right|\rightarrow0$.
The situation is simpler for the generalized eigenfunctions of $\hat{F}_{\nu}^{*}$
since they all belong to $\bigcap_{m>0}H^{m}=H^{\infty}=C^{\infty}\left(S^{1}\right)$. 
\item In the proof of Theorem \ref{th:gap_spectral}, we will obtain that
a general bound for $r_{s}\left(\hat{F}_{\nu}\right)$ (with no hypothesis
on $f$) is given by\begin{equation}
r_{s}\left(\hat{F}_{\nu}\right)\leq\frac{1}{\sqrt{E_{min}}}\exp\left(\frac{1}{2}\lim_{n\infty}\left(\frac{\log\mathcal{N}\left(n\right)}{n}\right)\right)+o\left(1\right)\label{eq:general_bound_rs}\end{equation}
where the function $\mathcal{N}\left(n\right)$ will be defined in
Eq.(\ref{eq:def_Nn}). This bound is similar to the bound given in
\cite[Theorem 1.1]{tsujii_05} by M. Tsujii.
\item In \cite[Theorem 1.2]{tsujii_05} M. Tsujii shows that \textbf{the
partially captive property}, i.e. $\lim_{n\infty}\left(\frac{\log\mathcal{N}\left(n\right)}{n}\right)=0$,
\textbf{is true for almost all functions} $\tau$.
\item From the definition of $\mathcal{N}\left(n\right)$ it is clear that
$\mathcal{N}\left(n\right)\leq k^{n}$ hence $\exp\left(\frac{1}{2}\lim_{n\infty}\left(\frac{\log\mathcal{N}\left(n\right)}{n}\right)\right)\leq\sqrt{k}$.
Also from the definition of $E_{min}$, it is clear that $E_{min}\leq k$
and therefore the upper bound in (\ref{eq:general_bound_rs}) is not
sharp since it does not give the obvious bound $r_{s}\left(\hat{F}_{\nu}\right)\leq1$
(see \cite[corollary 2]{fred-roy-sjostrand-07}). It is therefore
tempting to conjecture that for almost all functions $\tau$ Eq.(\ref{eq:upper_bound})
can be replaced by:\[
r_{s}\left(\hat{F}_{\nu}\right)\leq\frac{1}{\sqrt{k}}+o\left(1\right)\]

\item Notice that the above results say nothing about the existence of Ruelle
resonances $\lambda_{i}$. The work of F. Naud \cite{naud_08} are
the first results concerning  the existence of resonances $\lambda_{i}$. 
\item One observes numerically that for large $\nu\in\mathbb{R}$, the eigenvalues
$\lambda_{i}\left(\nu\right)$ repulse each other like eigenvalues
of random complex matrices. (See a movie on \url{http://www-fourier.ujf-grenoble.fr/~faure/articles}).
This suggests that many important questions of quantum chaos (e.g.
the conjecture of Random Matrices\cite{bohigas-89}) also concerns
the Ruelle resonances of partially hyperbolic dynamics in the semiclassical
limit. 
\item Remarks on numerical computation of the Ruelle resonances: one diagonalizes
the matrix which expresses the operator $\hat{F}_{\nu}$ in Fourier
basis $\varphi_{n}\left(x\right):=\exp\left(i2\pi nx\right)$, $n\in\mathbb{Z}$.
For the example of Figure \ref{fig:Spectrum-of-F} one gets  $\langle\varphi_{n'}|\hat{F}\varphi_{n}\rangle=e^{-i2\pi\frac{3}{4}\left(2n-n'\right)}J_{\left(2n-n'\right)}\left(\nu\right)$
where $J_{n}\left(x\right)$ is the Bessel function of first kind
\cite[9.1.21 p 360]{abramowitz}. Corollary 2 in \cite{fred-RP-06}
guaranties that the eigenvalues of the truncated matrix $\left|n\right|,\left|n'\right|\leq N$
converges towards the Ruelle resonances as $N\rightarrow\infty$.
\item One can proves \cite{arnoldi_09} that in the semi-classical limit
$\nu\rightarrow\infty$, the number of Ruelle resonances $\lambda_{i}$
(counting multiplicities) outside a fixed radius $\lambda$ is bounded
by a {}``\textbf{Weyl law}'':\[
\forall\lambda>0,\qquad\sharp\left\{ i\in\mathbb{N},\quad s.t.\,\,\left|\lambda_{i}\right|\geq\lambda\right\} \leq\left(\frac{\nu}{2\pi}\right)\mu\left(K\right)+o\left(\nu\right)\]
where $\mu\left(K\right)$ is the Lebesgue measure of the trapped
set $K$ defined later in Eq.(\ref{eq:def_K}). As usual in the semiclassical
theory of non selfadjoint operators, see \cite{sjostrand_90,sjostrand_07},
the Weyl law gives an upper bound for the density of resonances but
no lower bound. See discussions in \cite[section 3.1]{nonnenmacher_08}.
\end{itemize}

\subsection{Spectrum of $\hat{F}$ and dynamical correlation functions}

In this section, in order to give some {}``physical meaning'' to
the spectrum of $\hat{F}_{\nu}$, we recall relations between the
spectral results of Theorems \ref{th:spectrum_resonances},\ref{th:gap_spectral}
and the evolution of correlation functions \cite{baladi_livre_00}.
This will allow us to interpret the evolution and convergence of clouds
of points observed in Figure \ref{fig:Evolution}. 

Let $\nu\in\mathbb{Z}$. If $\psi_{1},\psi_{2}\in C^{\infty}\left(S^{1}\right)$,
the \textbf{correlation function at time} $n\in\mathbb{N}$ is defined
by:\[
C_{\psi_{2},\psi_{1}}\left(n\right):=\left(\hat{F}_{\nu}^{*n}\psi_{2},\psi_{1}\right)_{L^{2}}=\left(\psi_{2},\hat{F}_{\nu}^{n}\psi_{1}\right)_{L^{2}}\]
which represents the function $\psi_{2}$ evolved $n$ times by the
Perron-Frobenius operator $\hat{F}_{\nu}^{*}$ and tested against
the test function $\psi_{1}$. 

The first spectral result of Theorem \ref{th:spectrum_resonances}
implies that for any $\varepsilon>0$, and large $n$ (and assuming
that the eigenvalues $\left(\lambda_{i,\nu}\right)_{i}$ of $\hat{F}_{\nu}$
are simple for short; see \cite{fred-roy-sjostrand-07} for a more
extended discussion)\[
C_{\psi_{2},\psi_{1}}\left(n\right)=\sum_{\left|\lambda_{i,\nu}\right|>0}\lambda_{i,\nu}^{n}v_{i,\nu}\left(\overline{\psi}_{2}\right)w_{i,\nu}\left(\psi_{1}\right)+O_{\varepsilon}\left(\varepsilon^{n}\right)\]
If the conclusion of Theorem \ref{th:gap_spectral} holds, this implies
that for any $\rho$ such that $\frac{1}{\sqrt{E_{min}}}<\rho<1$,
there exists $\nu_{0}$ such that for any $\nu\geq\nu_{0}$, all the
eigenvalues of $\hat{F}_{\nu}$ are bounded: $\left|\lambda_{i,\nu}\right|<\rho<1,\forall i$.
This gives an exponential decay of correlations for $n\rightarrow\infty$
in these space $\hat{F}_{\nu}$:\[
C_{\psi_{2},\psi{}_{1}}\left(n\right)=\mathcal{O}\left(\rho^{n}\right)\]
It is known that if the function $\tau$ is not a co-boundary (i.e.
if the map $f$ is not equivalent to the trivial case $\tau=0$, as
explained in Appendix \ref{sub:Equivalence-classes}) then the map
$f$ is ergodic, which implies that all the eigenvalues $\lambda_{i,\nu}$
are strictly less than one: $\left|\lambda_{i,\nu}\right|<1,\forall\nu,\forall i$,
except for $\lambda_{0,0}=1$ associated to the eigenfunction $\varphi\left(x\right)=1$.
One deduces mixing property of the dynamics as observed in Figure
\ref{fig:Evolution}.

\section{\label{sec:Proof-of-theorem_spectral}Proof of theorem \ref{th:spectrum_resonances}
on resonances spectrum}

In this proof, we follow closely the proof of theorem 4 in \cite{fred-roy-sjostrand-07}
although we deal here with expanding map instead of hyperbolic map,
and this simplifies a lot, since we can work with ordinary Sobolev
spaces and not anisotropic Sobolev spaces. Here $\nu\in\mathbb{Z}$
is fixed.

\subsection{Dynamics on the cotangent space $T^{*}S^{1}$}

The first step is to realize that in order to study the spectral properties
of the transfer operator, we have to study the dynamics lifted on
the cotangent space. This basic idea has already been exploited in
\cite{fred-roy-sjostrand-07}. 

In Eq.(\ref{eq:def_E}), the map $E:S^{1}\rightarrow S^{1}$ is a
$k:1$ map, which means that every point $y\in S^{1}$ has $k$ inverses
denoted by $x_{\varepsilon}\in E^{-1}\left(y\right)$ and given explicitly
by \[
x_{\varepsilon}=E_{\varepsilon}^{-1}\left(y\right)=g^{-1}\left(\frac{y}{k}+\varepsilon\frac{1}{k}\right),\quad\mbox{with }\varepsilon=0,\ldots,k-1\]
We will denote the derivative by $E'\left(x\right):=dE/dx$.

\begin{figure}[h]
\begin{centering}
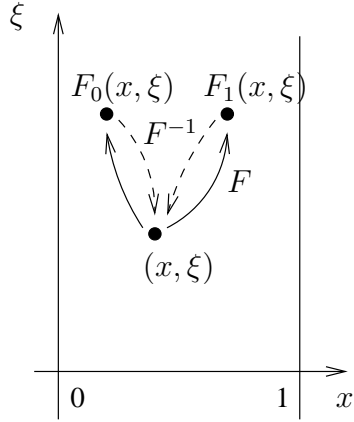
\par\end{centering}

\caption{\label{fig:The-map-F}This figure is for $k=2$. The map $F=\left\{ F_{0},\ldots,F_{k-1}\right\} $
is 1:k, and its inverse $F^{-1}$ is k:1 on $T^{*}S^{1}\equiv S^{1}\times\mathbb{R}$. }

\end{figure}

\vspace{0.cm}\begin{center}{\color{blue}\fbox{\color{black}\parbox{15cm}{
\begin{prop}\label{prop:F_OIF}

In Eq.(\ref{eq:def_F_hat}) $\hat{F}_{\nu}$ is a\textbf{ Fourier
integral operator} (FIO) acting on $C^{\infty}\left(S^{1}\right)$. 

The associated canonical transform on the cotangent space $\left(x,\xi\right)\in T^{*}S^{1}\equiv S^{1}\times\mathbb{R}$
is $k$-valued and given by:

\begin{equation}
F\left(x,\xi\right)=\left\{ F_{0}\left(x,\xi\right),\ldots,F_{k-1}\left(x,\xi\right)\right\} ,\qquad\left(x,\xi\right)\in S^{1}\times\mathbb{R}\label{eq:def_F}\end{equation}
where for any $\varepsilon=0,\ldots,k-1$,\begin{equation}
F_{\varepsilon}:\begin{cases}
x & \rightarrow x'_{\varepsilon}=E_{\varepsilon}^{-1}\left(x\right)=g^{-1}\left(\frac{1}{k}x+\varepsilon\frac{1}{k}\right)\\
\xi & \rightarrow\xi'_{\varepsilon}=E'\left(x'_{\varepsilon}\right)\xi=kg'\left(x'_{\varepsilon}\right)\xi\end{cases}\label{eq:expression_F}\end{equation}

Similarly the adjoint $\hat{F}^{*}$ is a FIO whose \textbf{canonical
transformation} is $F^{-1}$. See figure \ref{fig:The-map-F}.

\end{prop}
}}}\end{center}\vspace{0.cm}

The proof is just that the operator $\varphi\rightarrow\varphi\circ E$
on $C^{\infty}\left(S^{1}\right)$ is one of the simplest example
of Fourier integral operator, see \cite{martinez-01} example 2 p.150. 

The term $e^{i\nu\tau\left(x\right)}$ in Eq.(\ref{eq:def_F_hat})
does not contribute to the expression of $F$, since here $\nu$ is
considered as a fixed parameter, and therefore $e^{i\nu\tau\left(x\right)}$
acts as a pseudodifferential operator (equivalently as a FIO whose
canonical map is the identity).

The map $F$ is the map $E^{-1}$ lifted on the cotangent space $T^{*}S^{1}$
in the canonical way. Indeed, if we denote a point $\left(x,\xi\right)\in T^{*}S^{1}\equiv S^{1}\times\mathbb{R}$
then using the usual formula for differentials $y=E\left(x\right)=kg\left(x\right)$$\Rightarrow$$dy=E'\left(x\right)dx\Leftrightarrow\xi'=E'\left(x\right)\xi$,
we deduce the above expression for $F$.

\paragraph{Remarks}
\begin{itemize}
\item The \textbf{physical meaning} for $\hat{F}_{\varepsilon}$ being a
Fourier Integral Operator is that if $\varphi_{\left(x,\xi\right)}$
is a wave packet {}``micro-localized'' at position $\left(x,\xi\right)\in T^{*}S^{1}$
of phase space (this makes sense for $\xi\gg1$, and means that the
micro-support of $\varphi$ is $\left(x,\xi\right)$), then $\varphi':=\hat{F}_{\nu}\varphi_{\left(x,\xi\right)}$
will be a superposition of $k$ wave packets at positions $\left(x'_{\varepsilon},\xi'_{\varepsilon}\right)=F_{\varepsilon}\left(x,\xi\right)$,
$\varepsilon=0,\ldots,k-1$, i.e. with a very restricted micro-support,
controlled by the canonical map $F$.
\item Observe that the dynamics of the map $F$ on $S^{1}\times\mathbb{R}$
has a quite simple property: the zero section $\left\{ \left(x,\xi\right)\in S^{1}\times\mathbb{R},\,\xi=0\right\} $
is globally invariant and any other point with $\xi\neq0$ escapes
towards infinity $\left(\xi\rightarrow\pm\infty\right)$ in a controlled
manner:\begin{equation}
\left|\xi'_{\varepsilon}\right|\geq E_{min}\left|\xi\right|,\qquad\forall\varepsilon=0,\ldots,k-1\label{eq:expand}\end{equation}
where $E_{min}>1$ is given in (\ref{eq:Emin}).
\end{itemize}

\subsection{The escape function}

Let $m<0$ and define the $C^{\infty}$ function on $T^{*}S^{1}$:\[
A_{m}\left(x,\xi\right):=\left\langle \xi\right\rangle ^{m}\quad\in S^{m}\]
with $\left\langle \xi\right\rangle =\left(1+\xi^{2}\right)^{1/2}$.
$A_{m}$ decreases with $\left|\xi\right|$ and belongs to the symbol
class%
\footnote{See \cite{taylor_tome2} p.2. The \textbf{class of symbols} $S^{m}$,
with \textbf{order} $m\in\mathbb{R}$, consists of functions on the
cotangent space $A\in C^{\infty}\left(S^{1}\times\mathbb{R}\right)$
such that\[
\left|\partial_{\xi}^{\alpha}\partial_{x}^{\beta}A\right|_{\infty}\leq C_{\alpha,\beta}\left\langle \xi\right\rangle ^{m-\left|\alpha\right|},\qquad\left\langle \xi\right\rangle =\left(1+\xi^{2}\right)^{1/2}\]
} $S^{m}$.

Eq. (\ref{eq:expand}) implies that the function $A_{m}$ \textbf{decreases
strictly} along the trajectories of $F$ outside the zero section:\begin{equation}
\forall R>0,\forall\left|\xi\right|>R,\quad\forall\varepsilon=0,\ldots,k-1\qquad\frac{A_{m}\left(F_{\varepsilon}\left(x,\xi\right)\right)}{A_{m}\left(x,\xi\right)}\leq C^{\left|m\right|}<1,\quad\mbox{with }C=\sqrt{\frac{R^{2}+1}{R^{2}E_{min}+1}}<1\label{eq:A_decreases}\end{equation}

\begin{proof}
\[
\frac{A_{m}\left(F_{\varepsilon}\left(x,\xi\right)\right)}{A_{m}\left(x,\xi\right)}=\frac{\left(1+\xi^{2}\right)^{\left|m\right|/2}}{\left(1+\left(\xi'_{\varepsilon}\right)^{2}\right)^{\left|m\right|/2}}\leq\frac{\left(1+\xi^{2}\right)^{\left|m\right|/2}}{\left(1+E_{min}\xi^{2}\right)^{\left|m\right|/2}}\leq\left(\frac{1+R^{2}}{1+E_{min}R^{2}}\right)^{\left|m\right|/2}=C^{\left|m\right|}\]

\end{proof}
The symbol $A_{m}$ can be quantized into a pseudodifferential operator
$\hat{A}_{m}$ (PDO for short) which is self-adjoint and invertible
on $C^{\infty}\left(S^{1}\right)$ using the quantization rule (\cite{taylor_tome2}
p.2)\begin{equation}
\left(\hat{A}\varphi\right)\left(x\right)=\frac{1}{2\pi}\int A\left(x,\xi\right)e^{i\left(x-y\right)\xi}\varphi\left(y\right)dyd\xi,\label{eq:quantiz_rule}\end{equation}
but in our simple case, this is very explicit: in Fourier space, $\hat{A}_{m}$
is simply the multiplication by $\left\langle \xi\right\rangle ^{m}$. 

Remind that the \textbf{Sobolev space} $H^{m}\left(S^{1}\right)$
is defined by (\cite{taylor_tome1} p.271):\[
H^{m}\left(S^{1}\right):=\hat{A}_{m}^{-1}\left(L^{2}\left(S^{1}\right)\right)\]

The following commutative diagram\[
\begin{array}{ccc}
L^{2}\left(S^{1}\right) & \overset{\hat{Q}_{m}}{\rightarrow} & L^{2}\left(S^{1}\right)\\
\downarrow\hat{A}_{m}^{-1} & \circlearrowleft & \downarrow\hat{A}_{m}^{-1}\\
H^{m}\left(S^{1}\right) & \overset{\hat{F}_{\nu}}{\rightarrow} & H^{m}\left(S^{1}\right)\end{array}\]
 shows that $\hat{F}_{\nu}:H^{m}\left(S^{1}\right)\rightarrow H^{m}\left(S^{1}\right)$
is unitary equivalent to \[
\hat{Q}_{m}:=\hat{A}_{m}\hat{F}_{\nu}\hat{A}_{m}^{-1}\quad:L^{2}\left(S^{1}\right)\rightarrow L^{2}\left(S^{1}\right)\]
We will therefore study the operator $\hat{Q}_{m}$. Notice that $\hat{Q}_{m}$
is defined a priori on a dense domain ($C^{\infty}\left(S^{1}\right)$).
Define\begin{equation}
\hat{P}:=\hat{Q}_{m}^{*}\hat{Q}_{m}=\hat{A}_{m}^{-1}\left(\hat{F}_{\nu}^{*}\hat{A}_{m}^{2}\hat{F}_{\nu}\right)\hat{A}_{m}^{-1}=\hat{A}_{m}^{-1}\hat{B}\hat{A}_{m}^{-1}\label{eq:def_P}\end{equation}
where appears the operator \begin{equation}
\hat{B}:=\hat{F}_{\nu}^{*}\hat{A}_{m}^{2}\hat{F}_{\nu}\label{eq:def_B}\end{equation}

The Egorov Theorem will help us to treat this operator (see \cite{taylor_tome2}
p.24). This is a simple but crucial step in the proof: as explained
in \cite{fred-roy-sjostrand-07}, the Egorov theorem is the main Theorem
used in order to establish both the existence of a discrete spectrum
of resonances and properties of them. However there is a difference
with \cite{fred-roy-sjostrand-07}: for the expanding map we consider
here, the operator $\hat{F}_{\nu}$ is not invertible and the canonical
map $F$ is $k$-valued. Therefore we have to state the Egorov theorem
in an appropriate way (we restrict however the statement to our simple
context). 

\vspace{0.cm}\begin{center}{\color{blue}\fbox{\color{black}\parbox{15cm}{
\begin{lem}\textbf{\label{lem:Egorov-theorem}(Egorov theorem}). $\hat{B}:=\hat{F}_{\nu}^{*}\hat{A}_{m}^{2}\hat{F}_{\nu}$
is a pseudo-differential operator with symbol in $S^{m}$ given by:\begin{equation}
B\left(x,\xi\right)=\left(\sum_{\varepsilon=0,\ldots,k-1}\frac{1}{E'\left(x'_{\varepsilon}\right)}A_{m}^{2}\left(F_{\varepsilon}\left(x,\xi\right)\right)\right)+R\label{eq:Egorov}\end{equation}
with  $R\in S^{m-1}$ has a subleading order.

\end{lem}
}}}\end{center}\vspace{0.cm}

\begin{proof}

As we explained in Proposition \ref{prop:F_OIF}, $\hat{F}_{\nu}$
and $\hat{F}_{\nu}^{*}$ are Fourier integral operators (FIO) whose
canonical map are respectively $F$ and $F^{-1}$. The pseudodifferential
operator (PDO) $\hat{A}_{m}$ can also be considered as a FIO whose
canonical map is the identity. By composition we deduce that $\hat{B}=\hat{F}_{\nu}^{*}\hat{A}_{m}^{2}\hat{F}_{\nu}$
is a FIO whose canonical map is the identity since $F^{-1}\circ F=Id$.
See figure \ref{fig:The-map-F}. Therefore $\hat{B}$ is a PDO. Using
(\ref{eq:def_F_hat}), (\ref{eq:F_nu_*}) and (\ref{eq:expression_F})
we obtain that the principal symbol of $\hat{B}$ is \begin{equation}
\sum_{\varepsilon=0,\ldots,k-1}\frac{1}{E'\left(x'_{\varepsilon}\right)}A_{m}^{2}\left(F_{\varepsilon}\left(x,\xi\right)\right)\label{eq:ppcl_symbol_egorov}\end{equation}

\end{proof}

Remark: contrary to (\ref{eq:def_B}), $\hat{F}_{\nu}\hat{A}_{m}\hat{F}_{\nu}^{*}$
is not a PDO, but a FIO whose canonical map $F\circ F^{-1}$ is $k-$valued
(see figure \ref{fig:The-map-F}).

Now by\textbf{ theorem of composition of PDO} (\cite{taylor_tome2}
p.11), (\ref{eq:def_P}) and (\ref{eq:Egorov}) imply that $\hat{P}$
is a PDO of order $0$ with principal symbol:\[
P\left(x,\xi\right)=\frac{B\left(x,\xi\right)}{A_{m}^{2}\left(x,\xi\right)}=\left(\sum_{\varepsilon=0,\ldots,k-1}\frac{1}{E'\left(x'_{\varepsilon}\right)}\frac{A_{m}^{2}\left(F_{\varepsilon}\left(x,\xi\right)\right)}{A_{m}^{2}\left(x,\xi\right)}\right)\]
The estimate (\ref{eq:A_decreases}) together with (\ref{eq:Emin})
give the following upper bound\[
\forall\left|\xi\right|>R,\quad\left|P\left(x,\xi\right)\right|\leq C^{2\left|m\right|}\sum_{\varepsilon=0,\ldots,k-1}\frac{1}{E'\left(x'_{\varepsilon}\right)}\leq C^{2\left|m\right|}\frac{k}{E_{min}}\]
(This upper bound goes to zero as $m\rightarrow-\infty$). From\textbf{
$L^{2}$-continuity theorem for PDO }we deduce that for any $\alpha>0$
(see \cite{fred-roy-sjostrand-07} Lemma 38)\[
\hat{P}=\hat{k}_{\alpha}+\hat{p}_{\alpha}\]
with $\hat{k}_{\alpha}$ a smoothing operator (hence compact) and
$\left\Vert \hat{p}_{\alpha}\right\Vert \leq C^{2\left|m\right|}\frac{k}{E_{min}}+\alpha$.
If $\hat{Q}_{m}=\hat{U}\left|\hat{Q}\right|$ is the polar decomposition
of $\hat{Q}_{m}$, with $\hat{U}$ unitary, then from (\ref{eq:def_P})
$\hat{P}=\left|\hat{Q}\right|^{2}\Leftrightarrow\left|\hat{Q}\right|=\sqrt{\hat{P}}$
and the spectral theorem (\cite{taylor_tome2} p.75) gives that $\left|\hat{Q}\right|$
has a similar decomposition\[
\left|\hat{Q}\right|=\hat{k}'_{\alpha}+\hat{q}_{\alpha}\]
with $\hat{k}'_{\alpha}$ smoothing and $\left\Vert \hat{q}_{\alpha}\right\Vert \leq C^{\left|m\right|}\sqrt{\frac{k}{E_{min}}}+\alpha$,
with any $\alpha>0$. Since $\left\Vert \hat{U}\right\Vert =1$ we
deduce a similar decomposition for $\hat{Q}_{m}=\hat{U}\left|\hat{Q}\right|:L^{2}\left(S^{1}\right)\rightarrow L^{2}\left(S^{1}\right)$
and we deduce (\ref{eq:decomp_F_nu}) and (\ref{eq:bound_R}) for
$\hat{F}_{\nu}:H^{m}\rightarrow H^{m}$. We also use the fact that
$C\rightarrow1/E_{min}$ for $R\rightarrow\infty$ in (\ref{eq:A_decreases}). 

The fact that the eigenvalues $\lambda_{i}$ and their generalized
eigenspaces do not depend on the choice of space $H^{m}$ is due to
density of Sobolev spaces. We refer to the argument given in the proof
of corollary 1 in \cite{fred-roy-sjostrand-07}. This finishes the
proof of Theorem \ref{th:spectrum_resonances}.

\section{Proof of theorem \ref{th:gap_spectral} on spectral gap}

We will follow steps by steps the same analysis as in the previous
section. The main difference now is that in Theorem \ref{th:gap_spectral},
$\nu\gg1$ is a semi-classical parameter. In other words, we just
perform a linear rescaling in cotangeant space: $\xi_{h}:=\hbar\xi$
with \[
\hbar:=\frac{1}{\nu}\ll1.\]
Therefore, our quantization rule for a symbol $A\left(x,\xi_{h}\right)$,
Eq.(\ref{eq:quantiz_rule}) writes now (see \cite{martinez-01} p.22)\begin{equation}
\left(\hat{A}\varphi\right)\left(x\right)=\frac{1}{2\pi\hbar}\int A\left(x,\xi\right)e^{i\left(x-y\right)\xi_{h}/\hbar}\varphi\left(y\right)dyd\xi_{h}\label{eq:quantiz_rule_n}\end{equation}
For simplicity we will write $\xi$ for $\xi_{h}$ below.

\subsection{Dynamics on the cotangent space $T^{*}S^{1}$}

In Eq.(\ref{eq:def_F_hat}) the multiplicative term $e^{i\nu\tau\left(x\right)}=e^{i\tau\left(x\right)/\hbar}$
acts now as a Fourier integral operator (FIO) and contributes to the
transport (it was merely a P.D.O. for theorem \ref{th:spectrum_resonances}
in Section \ref{sec:Proof-of-theorem_spectral} when $\nu$ was fixed).
Its associated canonical transformation on $T^{*}S^{1}=S^{1}\times\mathbb{R}$
is $\left(x,\xi\right)\rightarrow\left(x,\xi+\frac{d\tau}{dx}\left(x\right)\right)$
(this is a direct consequence of stationary phase approximation in
Fourier transform see \cite[Examples 1,2 p.150]{martinez-01}). We
obtain:

\vspace{0.cm}\begin{center}{\color{blue}\fbox{\color{black}\parbox{15cm}{
\begin{prop}

In Eq.(\ref{eq:def_F_hat}) $\hat{F}_{\nu}$ is a \textbf{semi-classical
Fourier integral operator} acting on $C^{\infty}\left(S^{1}\right)$(with
semi-classical parameter $\hbar:=1/\nu\ll1$). The associated \textbf{canonical
transformation} on the cotangent space $\left(x,\xi\right)\in T^{*}S^{1}\equiv S^{1}\times\mathbb{R}$
is $k$-valued and given by:

\begin{equation}
F\left(x,\xi\right)=\left\{ F_{0}\left(x,\xi\right),\ldots,F_{k-1}\left(x,\xi\right)\right\} ,\qquad\left(x,\xi\right)\in S^{1}\times\mathbb{R}\label{eq:def_F_nu}\end{equation}
\begin{equation}
F_{\varepsilon}:\begin{cases}
x & \rightarrow x'_{\varepsilon}=E^{-1}\left(x\right)\\
\xi & \rightarrow\xi'_{\varepsilon}=E'\left(x'_{\varepsilon}\right)\xi+\frac{d\tau}{dx}\left(x'_{\varepsilon}\right)\end{cases},\qquad\varepsilon=0,\ldots,k-1\label{eq:expression_F_nu}\end{equation}

Similarly $\hat{F}^{*}$ is a FIO whose canonical transformation is
$F^{-1}$. 

\end{prop}
}}}\end{center}\vspace{0.cm}

Notice that for simplicity we have kept the same notation for the
canonical transformation $F$ although it differs from (\ref{eq:expression_F}).

Since the map $F$ is $k$-valued, a trajectory is a tree. Let us
precise the notation: 

\vspace{0.cm}\begin{center}{\color{red}\fbox{\color{black}\parbox{15cm}{
\begin{Def}

For $\varepsilon=\left(\ldots\varepsilon_{3},\varepsilon_{2},\varepsilon_{1}\right)\in\left\{ 0,\ldots,k-1\right\} ^{\mathbb{N}^{*}}$,
a point $\left(x,\xi\right)\in S^{1}\times\mathbb{R}$ and time $n\in\mathbb{N}^{*}$
let us denote:

\begin{equation}
F_{\varepsilon}^{n}\left(x,\xi\right):=F_{\varepsilon_{n}}F_{\varepsilon_{n-1}}\ldots F_{\varepsilon_{1}}\left(x,\xi\right)\label{def:trajectory}\end{equation}

For a given sequence $\varepsilon\in\left\{ 0,\ldots,k-1\right\} ^{\mathbb{N}^{*}}$,
a \textbf{trajectory }issued from the point $\left(x,\xi\right)$
is $\left\{ F_{\varepsilon}^{n}\left(x,\xi\right),\,\, n\in\mathbb{N}\right\} $.

\end{Def}
}}}\end{center}\vspace{0.cm}

Notice that at time $n\in\mathbb{N}$, there are $k^{n}$ points issued
from a given point $\left(x,\xi\right)$:\begin{equation}
F^{n}\left(x,\xi\right):=\left\{ F_{\varepsilon}^{n}\left(x,\xi\right),\quad\varepsilon\in\left\{ 0,\ldots,k-1\right\} ^{n}\right\} \end{equation}

The new term $\frac{d\tau}{dx}\left(x'_{\varepsilon}\right)$ in the
expression of $\xi'_{\varepsilon}$, Eq.(\ref{eq:expression_F_nu}),
complicates significantly the dynamics near the zero section $\xi=0$.
However a trajectory from an initial point with $\left|\xi\right|$
large enough still escape towards infinity:

\vspace{0.cm}\begin{center}{\color{blue}\fbox{\color{black}\parbox{15cm}{
\begin{lem}

For any $1<\kappa<E_{min}$, there exists $R\geq0$ such that for
any $\left|\xi\right|>R$, any $\varepsilon=0,\ldots k-1$, \begin{equation}
\left|\xi'_{\varepsilon}\right|>\kappa\left|\xi\right|\label{eq:expand_nu}\end{equation}

\end{lem}
}}}\end{center}\vspace{0.cm}

\begin{proof}
\small 

From (\ref{eq:expression_F_nu}), one has $\xi'_{\varepsilon}=E'\left(x'_{\varepsilon}\right)\xi+\tau'\left(x'_{\varepsilon}\right)$,
so $\xi'_{\varepsilon}-\kappa\xi=\left(E'\left(x'_{\varepsilon}\right)-\kappa\right)\xi+\tau'\left(x'_{\varepsilon}\right)\geq\left(E_{min}-\kappa\right)\xi+\min\tau'>0$
if $\xi>-\frac{\min\tau'}{\left(E_{min}-\kappa\right)}\geq0$, and
similarly $\xi'_{\varepsilon}-\kappa\xi\leq\left(E_{min}-\kappa\right)\xi+\max\tau'<0$
if $\xi<-\frac{\max\tau'}{\left(E_{min}-\kappa\right)}$.\end{proof}
\normalsize

We will denote the set:\begin{equation}
\boxed{\mathcal{Z}:=S^{1}\times\left[-R,R\right]}\label{eq:def_Z}\end{equation}
outside of which trajectories escape in a controlled manner (\ref{eq:expand_nu}).
See figure \ref{fig:The-trapped-set}.

\subsection{The trapped set $K$}

We will be interested now in the trajectories of $F$ which do not
escape towards infinity.

\begin{figure}[h]
\begin{centering}
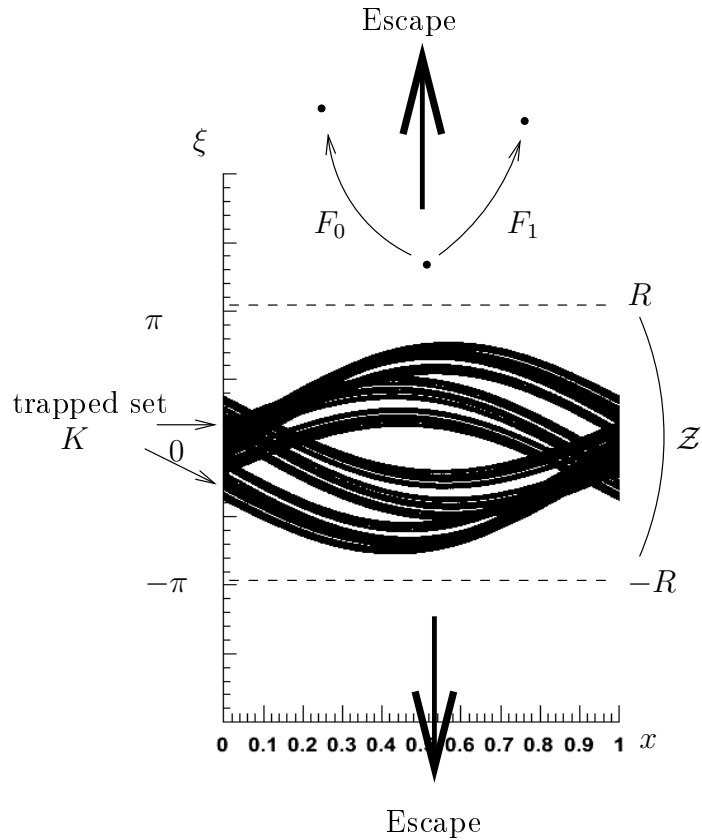
\par\end{centering}

\caption{\label{fig:The-trapped-set}The trapped set $K$ in the cotangent
space $S^{1}\times\mathbb{R}$. We have chosen here $E\left(x\right)=2x$
and $\tau\left(x\right)=\cos\left(2\pi x\right)$.}

\end{figure}

\vspace{0.cm}\begin{center}{\color{red}\fbox{\color{black}\parbox{15cm}{
\begin{Def}

We define\underbar{ }\textbf{the trapped set}

\begin{equation}
\boxed{K:=\bigcap_{n\in\mathbb{N}}\left(F^{-1}\right)^{n}\left(\mathcal{Z}\right)}\label{eq:def_K}\end{equation}
which contains points for which a trajectory at least does not escape
towards infinity. See figure \ref{fig:The-trapped-set}. The definition
of $K$ does not depend on the compact set $\mathcal{Z}$ (if $\mathcal{Z}$
is chosen large enough).

\end{Def}
}}}\end{center}\vspace{0.cm}

Since the map $F$ is multivalued, some trajectories may escape from
the trapped set. We will need a characterization of how many such
trajectories succeed to escape:

For $n\in\mathbb{N}$, let\begin{equation}
\boxed{\mathcal{N}\left(n\right):=\max_{\left(x,\xi\right)}\sharp\left\{ F_{\varepsilon}^{n}\left(x,\xi\right)\in\mathcal{Z},\quad\varepsilon\in\left\{ 0,\ldots,k-1\right\} ^{n}\right\} }\label{eq:def_Nn}\end{equation}

\begin{figure}[h]
\begin{centering}
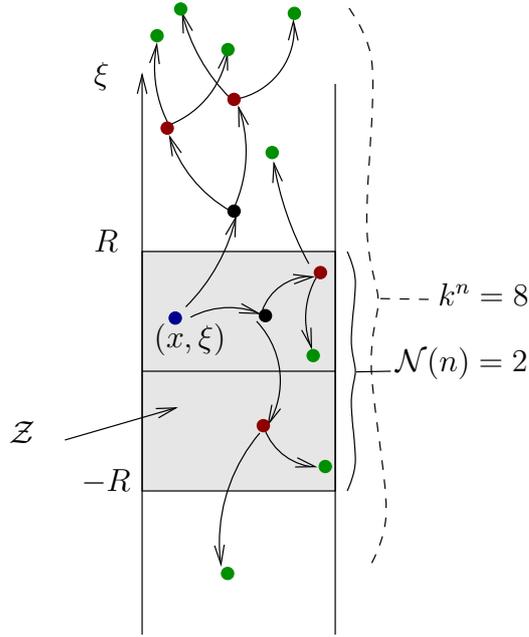
\par\end{centering}

\caption{\label{fig:partially_captive}This Figure illustrates the \textbf{trajectories}
$F_{\varepsilon}^{n}\left(x,\xi\right)$ issued from an initial point
$\left(x,\xi\right)$. Here $k=2$ and $n=3$. The property for the
map $F$ of being {}``partially captive'' according to definition
\ref{def:partially_captive} is related to the number of points $\mathcal{N}\left(n\right)$
which do not escape from the compact zone $\mathcal{Z}$ after time
$n$.}

\end{figure}

See Figure \ref{fig:partially_captive} for an illustration of $\mathcal{N}\left(n\right)$.
Of course \emph{$\mathcal{N}\left(n\right)\leq k^{n}$.} 

\vspace{0.cm}\begin{center}{\color{red}\fbox{\color{black}\parbox{15cm}{
\begin{Def}

\label{def:partially_captive}

The map $F$ (or $f$) is \textbf{partially captive} if\begin{equation}
\frac{\log\mathcal{N}\left(n\right)}{n}\underset{n\rightarrow\infty}{\longrightarrow}0\label{eq:def_partially_captive}\end{equation}
This property is the hypothesis of Theorem \ref{th:gap_spectral}.

\end{Def}
}}}\end{center}\vspace{0.cm}

\paragraph{Remarks}
\begin{itemize}
\item {}``$F$ partially captive'' means that most of the trajectories
escape from the trapped set $K$. See figure \ref{fig:partially_captive}.
Another description of the trapped set $K$ and of the partially captive
property will be given in Appendix \ref{sec:Description-of-K}. Notice
that the function $\mathcal{N}\left(n\right)$, Eq.(\ref{eq:def_Nn})
depends on the set $\mathcal{Z}$ but property (\ref{eq:def_partially_captive})
does not.
\item If the function $\tau$ is trivial in (\ref{eq:expression_F_nu}),
i.e. $\tau=0$ , then obviously all the trajectories issued from a
point $\left(x,\xi\right)$ on the line $\xi=0$ remains on this line
(the trapped set). Therefore \[
\sharp\left\{ F_{\varepsilon}^{n}\left(x,\xi\right)\in\mathcal{Z},\varepsilon\in\left\{ 0,1,\ldots,k-1\right\} ^{n}\right\} =k^{n}\]
 and the map $F$ is not partially captive (but could be called {}``totally
captive''). This is also true if the function $\tau$ is a {}``co-boundary'',
i.e. if $\tau\left(x\right)=\eta\left(E\left(x\right)\right)-\eta\left(x\right)$
with $\eta\in C^{\infty}\left(S^{1}\right)$ as discussed in Appendix
\ref{sub:Equivalence-classes}.
\item M. Tsujii has studied a dynamical system very similar to (\ref{eq:expression_F_nu})
in \cite{tsujii_01}, but this model is not volume preserving. He
establishes there that the SRB measure on the trapped set is absolutely
continuous for almost every $\tau$.
\end{itemize}

\subsection{The escape function}

Let $m<0$ and consider the $C^{\infty}$ function on $T^{*}S^{1}$:\begin{eqnarray*}
A_{m}\left(x,\xi\right) & := & \left\langle \xi\right\rangle ^{m}\qquad\mbox{for }\left|\xi\right|>R+\eta\\
 & := & 1\qquad\mbox{for }\xi\leq R\end{eqnarray*}
where $\eta>0$ is small and with $\left\langle \xi\right\rangle :=\left(1+\xi^{2}\right)^{1/2}$.
$A_{m}$ decreases with $\left|\xi\right|$ and belongs to the symbol
class $S^{m}$.

Eq. (\ref{eq:expand_nu}) implies that the function $A_{m}$ \textbf{decreases
strictly} along the trajectories of $F$ outside the trapped set (similarly
to Eq.(\ref{eq:A_decreases})):\begin{equation}
\forall\left|\xi\right|>R,\quad\forall\varepsilon=0,\ldots,k-1\qquad\frac{A_{m}\left(F_{\varepsilon}\left(x,\xi\right)\right)}{A_{m}\left(x,\xi\right)}\leq C^{\left|m\right|}<1,\quad\mbox{with }C=\sqrt{\frac{R^{2}+1}{\kappa R^{2}+1}}<1\label{eq:A_decreases_2}\end{equation}

And for any point we have the general bound:\begin{equation}
\forall\left(x,\xi\right)\in T^{*}S^{1},\qquad\frac{A_{m}\left(F_{\varepsilon}\left(x,\xi\right)\right)}{A_{m}\left(x,\xi\right)}\leq1.\label{eq:general_bound_A}\end{equation}
Using the quantization rule (\ref{eq:quantiz_rule_n}), the symbol
$A_{m}$ can be quantized giving a pseudodifferential operator $\hat{A}_{m}$
which is self-adjoint and invertible on $C^{\infty}\left(S^{1}\right)$.
In our case $\hat{A}_{m}$ is simply a multiplication operator by
$A_{m}\left(\xi\right)$ in Fourier space. 

Let us consider the (usual) Sobolev space \[
H^{m}\left(S^{1}\right):=\hat{A}_{m}^{-1}\left(L^{2}\left(S^{1}\right)\right)\]
Then $\hat{F}_{\nu}:H^{m}\left(S^{1}\right)\rightarrow H^{m}\left(S^{1}\right)$
is unitary equivalent to \[
\hat{Q}:=\hat{A}_{m}\hat{F}_{\nu}\hat{A}_{m}^{-1}\quad:L^{2}\left(S^{1}\right)\rightarrow L^{2}\left(S^{1}\right)\]
Let $n\in\mathbb{N}^{*}$ (a fixed time which will be made large at
the end of the proof) and define

\begin{equation}
\hat{P}^{\left(n\right)}:=\hat{Q}^{*n}\hat{Q}^{n}=\hat{A}_{m}^{-1}\hat{F}_{\nu}^{*n}\hat{A}_{m}^{2}\hat{F}_{\nu}^{n}\hat{A}_{m}^{-1}\label{eq:def_Pn}\end{equation}
Using \textbf{Egorov theorem} (the semi-classical version of\textbf{
}Lemma \ref{lem:Egorov-theorem}) and Theorem of composition of PDO,
we obtain that $\hat{P}^{\left(n\right)}$ is a PDO of order $0$
with principal symbol\begin{equation}
P^{\left(n\right)}\left(x,\xi\right)=\left(\sum_{\varepsilon\in\left\{ 0,\ldots,k-1\right\} ^{n}}\frac{1}{E'_{n}\left(x\right)}\frac{A_{m}^{2}\left(F_{\varepsilon}^{n}\left(x,\xi\right)\right)}{A_{m}^{2}\left(x,\xi\right)}\right)\label{eq:e1}\end{equation}
where $E'_{n}\left(x\right):=\prod_{j=1}^{n}E'\left(E_{\varepsilon_{j}}^{-j}\left(x\right)\right)$
is the expanding rate of the trajectory at time $n$. Eq.(\ref{eq:Emin})
implies that $E'_{n}\left(x\right)\geq E_{min}^{n}$. Now we will
bound this (positive) symbol from above, considering different cases
for the trajectory $F_{\varepsilon}^{n}\left(x,\xi\right)$, as illustrated
on Figure \ref{fig:partially_captive}. 
\begin{enumerate}
\item If $\left(x,\xi\right)\notin\mathcal{Z}$ then (\ref{eq:A_decreases_2})
gives \begin{equation}
\frac{A^{2}\left(F_{\varepsilon}^{n}\left(x,\xi\right)\right)}{A^{2}\left(x,\xi\right)}=\frac{A^{2}\left(F_{\varepsilon}^{n}\left(x,\xi\right)\right)}{A^{2}\left(F_{\varepsilon}^{n-1}\left(x,\xi\right)\right)}\ldots\frac{A^{2}\left(F_{\varepsilon}\left(x,\xi\right)\right)}{A^{2}\left(x,\xi\right)}\leq\left(C^{2\left|m\right|}\right)^{n}\label{eq:bound1}\end{equation}
therefore \[
P^{\left(n\right)}\left(x,\xi\right)\leq\frac{k^{n}}{E_{min}^{n}}\left(C^{2\left|m\right|}\right)^{n}\]

\item If $\left(x,\xi\right)\in\mathcal{Z}$ but $F_{\varepsilon}^{n-1}\left(x,\xi\right)\notin\mathcal{Z}$
then $\frac{\left(A^{2}\circ F_{\varepsilon}^{n}\right)\left(x,\xi\right)}{\left(A^{2}\circ F_{\varepsilon}^{n-1}\right)\left(x,\xi\right)}\leq C^{2\left|m\right|}$
from (\ref{eq:A_decreases_2}). Using also (\ref{eq:general_bound_A})
we have \begin{equation}
\frac{A^{2}\left(F_{\varepsilon}^{n}\left(x,\xi\right)\right)}{A^{2}\left(x,\xi\right)}=\frac{A^{2}\left(F_{\varepsilon}^{n}\left(x,\xi\right)\right)}{A^{2}\left(F_{\varepsilon}^{n-1}\left(x,\xi\right)\right)}\ldots\frac{A^{2}\left(F_{\varepsilon}\left(x,\xi\right)\right)}{A^{2}\left(x,\xi\right)}\leq C^{2\left|m\right|}\label{eq:e2}\end{equation}

\item In the other cases ($\left(x,\xi\right)\in\mathcal{Z}$ and $F_{\varepsilon}^{n-1}\left(x,\xi\right)\in\mathcal{Z}$)
we can only use (\ref{eq:general_bound_A}) to bound: \begin{equation}
\frac{A^{2}\left(F_{\varepsilon}^{n}\left(x,\xi\right)\right)}{A^{2}\left(x,\xi\right)}\leq1\label{eq:e3}\end{equation}

\end{enumerate}
From definition (\ref{eq:def_Nn}) we have

\[
\sharp\left\{ F_{\varepsilon}^{n-1}\left(x,\xi\right)\in\mathcal{Z},\qquad\varepsilon\in\left\{ 0,1\right\} ^{n}\right\} \leq\mathcal{N}\left(n-1\right).\]
For $\left(x,\xi\right)\in\mathcal{Z}$, we split the sum Eq.(\ref{eq:e1})
accordingly to cases 1,2 or 3 above. Notice that $\left(C^{2\left|m\right|}\right)^{n}\leq C^{2\left|m\right|}$.
This gives \begin{equation}
P^{\left(n\right)}\left(x,\xi\right)\leq\frac{1}{E_{min}^{n}}\left(\left(k^{n}-\mathcal{N}\left(n-1\right)\right)C^{2\left|m\right|}+\mathcal{N}\left(n-1\right)\right)\leq\mathcal{B}\label{eq:bound2}\end{equation}

with the bound\[
\mathcal{B}:=\left(\frac{k}{E_{min}}\right)^{n}C^{2\left|m\right|}+\frac{\mathcal{N}\left(n-1\right)}{E_{min}^{n}}\]
Then\[
\limsup_{\left(x,\xi\right)}\left|P^{\left(n\right)}\left(x,\xi\right)\right|\leq\mathcal{B}\]
With \textbf{$L^{2}$continuity theorem} for pseudodifferential operators
this implies that in the limit $\hbar\rightarrow0$\begin{equation}
\left\Vert \hat{P}^{\left(n\right)}\right\Vert \leq\mathcal{B}+\mathcal{O}_{n}\left(\hbar\right)\end{equation}

Polar decomposition of $\hat{Q}^{n}$ gives\[
\left\Vert \hat{Q}^{n}\right\Vert \leq\left\Vert \left|\hat{Q}^{n}\right|\right\Vert =\sqrt{\left\Vert \hat{P}^{\left(n\right)}\right\Vert }\leq\left(\mathcal{B}+\mathcal{O}_{n}\left(\hbar\right)\right)^{1/2}\]
Then for any $n$ the spectral radius of $\hat{Q}$ satisfies \cite[p.192]{reed-simon1}
\[
r_{s}\left(\hat{Q}\right)\leq\left\Vert \hat{Q}^{n}\right\Vert ^{1/n}\leq\left(\mathcal{B}+\mathcal{O}_{n}\left(\hbar\right)\right)^{1/2n}\]
 Also notice that\[
\left(\frac{\mathcal{N}\left(n-1\right)}{E_{min}^{n}}\right)^{1/2n}=\frac{1}{\sqrt{E_{min}}}\exp\left(\frac{1}{2n}\log\mathcal{N}\left(n-1\right)\right)\]

We let $\hbar\rightarrow0$ first, and after $m\rightarrow-\infty$
giving $C^{\left|m\right|}\rightarrow0$, and also we let $n\rightarrow\infty$.
Then for $\hbar=1/\nu\rightarrow0$ we have%
\footnote{It can be shown that $\log\mathcal{N}\left(n\right)$ is sub-additive
and therefore $\lim_{n\infty}\mbox{inf}\left(\frac{\log\mathcal{N}\left(n\right)}{n}\right)=\lim_{n\infty}\left(\frac{\log\mathcal{N}\left(n\right)}{n}\right)$,\cite[p.217, ex.11]{reed-simon1}%
}:\begin{equation}
r_{s}\left(\hat{Q}\right)\leq\sqrt{\frac{1}{E_{min}}\exp\left(\lim_{n\infty}\mbox{inf}\left(\frac{\log\mathcal{N}\left(n\right)}{n}\right)\right)}+o\left(1\right).\label{eq:bound_rs}\end{equation}
\[
\]
If we make the assumption that $F$ be \textbf{partially captive},
Eq.(\ref{eq:def_partially_captive}), we get that for $\hbar=1/\nu\rightarrow0$,\[
r_{s}\left(\hat{Q}\right)\leq\frac{1}{\sqrt{E_{min}}}+o\left(1\right).\]
We have finished the proof of Theorem \ref{th:gap_spectral}.

\appendix
\newpage{}

\section{\label{sub:Equivalence-classes}Equivalence classes of dynamics}

Let us make a simple and well known observation about equivalent classes
of dynamics. The map $f$ we consider in Eq.(\ref{eq:def_f}) depends
on $k\in\mathbb{N}$ and on the functions $E:S^{1}\rightarrow S^{1}$,
$\tau:S^{1}\rightarrow\mathbb{R}$. To emphasize this dependence,
we denote $f_{\left(E,\tau\right)}$. The transfer operator (\ref{eq:def_F_hat})
is also denoted by $\hat{F}_{\left(E,\tau\right)}$.

In this Appendix we characterize an equivalence class of functions
$\left(E,\tau\right)$ such that in a given equivalence class the
maps $f_{\left(E,\tau\right)}$ are $C^{\infty}$ conjugated together,
the transfer operators $\hat{F}_{\left(E,\tau\right)}$ are also conjugated
and the resonances spectrum are therefore the same. 

Let $\eta:S^{1}\rightarrow\mathbb{R}$ be a smooth function. Let us
consider the map $T:S^{1}\times S^{1}\rightarrow S^{1}\times S^{1}$
defined by\[
T\left(x,s\right)=\left(x,s+\frac{1}{2\pi}\eta\left(x\right)\right)\]
Then using (\ref{eq:def_f}) one gets that:\[
\left(T^{-1}\circ f_{\left(E,\tau\right)}\circ T\right)\left(x,s\right)=\left(E\left(x\right),s+\frac{1}{2\pi}\left(\tau\left(x\right)+\eta\left(x\right)-\eta\left(E\left(x\right)\right)\right)\right)\]
Therefore\[
\boxed{\left(T^{-1}\circ f_{\left(E,\tau\right)}\circ T\right)=f_{\left(E,\zeta\right)}}\]
i.e. $f_{\left(E,\zeta\right)}\sim f_{\left(E,\tau\right)}$, with\[
\boxed{\zeta=\tau+\left(\eta-\eta\circ E\right).}\]
The function $\tau$ has been modified by a {}``\textbf{co-boundary
term}'' (\cite{katok_hasselblatt}, p.100). 

With (\ref{eq:def_F_hat}) we also obtain that the transfer operator
$\hat{F}_{\left(E,\zeta\right)}$ of $f_{\left(E,\zeta\right)}$ on
$C^{\infty}\left(S^{1}\right)$ is given by\begin{equation}
\boxed{\hat{F}_{\left(E,\zeta\right)}=\hat{\chi}\hat{F}_{\left(E,\tau\right)}\hat{\chi}^{-1}}\label{eq:conjugated_op}\end{equation}
with the operator $\hat{\chi}:C^{\infty}\left(S^{1}\right)\rightarrow C^{\infty}\left(S^{1}\right)$
defined by: \[
\boxed{\left(\hat{\chi}\varphi\right)\left(x\right)=e^{i\nu\eta\left(x\right)}\varphi\left(x\right).}\]

\begin{proof}
$\left(\hat{\chi}\hat{F}_{\left(E,\tau\right)}\hat{\chi}^{-1}\varphi\right)\left(x\right)=\left(\varphi\left(E\left(x\right)\right)e^{-i\nu\eta\left(E\left(x\right)\right)}\right)e^{i\nu\tau\left(x\right)}e^{i\nu\eta\left(x\right)}=\left(\hat{F}_{\left(E,\zeta\right)}\varphi\right)\left(x\right)$.
\end{proof}
The conjugation (\ref{eq:conjugated_op}) immediately implies that
$\hat{F}_{\left(E,\zeta\right)}$ and $\hat{F}_{\left(E,\tau\right)}$
have the same spectrum of Ruelle resonances.

Observe that $\hat{\chi}$ is a O.I.F whose associated canonical transformation
on $T^{*}S^{1}\equiv S^{1}\times\mathbb{R}$ is given by ($\nu\gg1$
is considered as a semi-classical parameter):\[
\chi:\left(x,\xi\right)\in\left(S^{1}\times\mathbb{R}\right)\rightarrow\left(x,\xi+\frac{d\eta}{dx}\right)\in\left(S^{1}\times\mathbb{R}\right).\]
Therefore at the level of canonical transforms on $T^{*}S^{1}$:\begin{equation}
\boxed{F_{\left(E,\zeta\right)}=\chi\circ F_{\left(E,\tau\right)}\circ\chi^{-1}}\label{eq:conjug_F}\end{equation}

The conjugation (\ref{eq:conjug_F}) implies in particular that the
corresponding trapped sets (\ref{eq:def_K}) are related by \[
\boxed{K_{\left(E,\zeta\right)}=\chi\left(K_{\left(E,\tau\right)}\right)}\]

\section{\label{sec:Description-of-K}Description of the trapped set}

In this section we provide further description of the trapped set
$K$ defined in Eq.(\ref{eq:def_K}) as well as the dynamics of the
canonical map $F$ restricted on it.

\subsection{Dynamics on the cover $\mathbb{R}^{2}$}

The dynamics of $F$ on the cylinder $T^{*}S^{1}=S^{1}\times\mathbb{R}$
has been defined in Eq.(\ref{eq:expression_F_nu}). This is a multivalued
map. It is convenient to consider the \textbf{lifted dynamics on the
cover} $\mathbb{R}^{2}$ which is a diffeomorphism given by\begin{equation}
\tilde{F}:\begin{cases}
x & \rightarrow x'=E^{-1}\left(x\right)\\
\xi & \rightarrow\xi'=E'\left(x'\right)\xi+\frac{d\tau}{dx}\left(x'\right)\end{cases}\label{eq:F_lifted}\end{equation}

where $E=kg:\mathbb{R}\rightarrow\mathbb{R}$ is the map (\ref{eq:def_E})
lifted on $\mathbb{R}$. It is invertible from (\ref{eq:Emin}). Let
us suppose for simplicity that $E\left(0\right)=0$.

\begin{figure}[h]
\begin{centering}
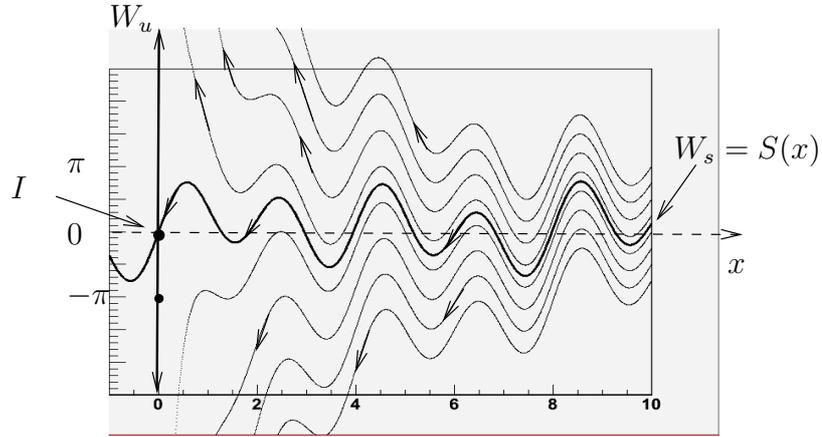
\par\end{centering}

\caption{\label{fig:separatrice}The fixed point $I=\left(0,\xi_{0}\right)$,
the stable manifold $W_{s}$ and unstable manifold $W_{u}$ of the
lifted map $\tilde{F}$, Eq.(\ref{eq:F_lifted}), in the example $E\left(x\right)=2x$,
$\tau\left(x\right)=\cos\left(2\pi x\right)$.}

\end{figure}

\begin{figure}[h]
\begin{centering}
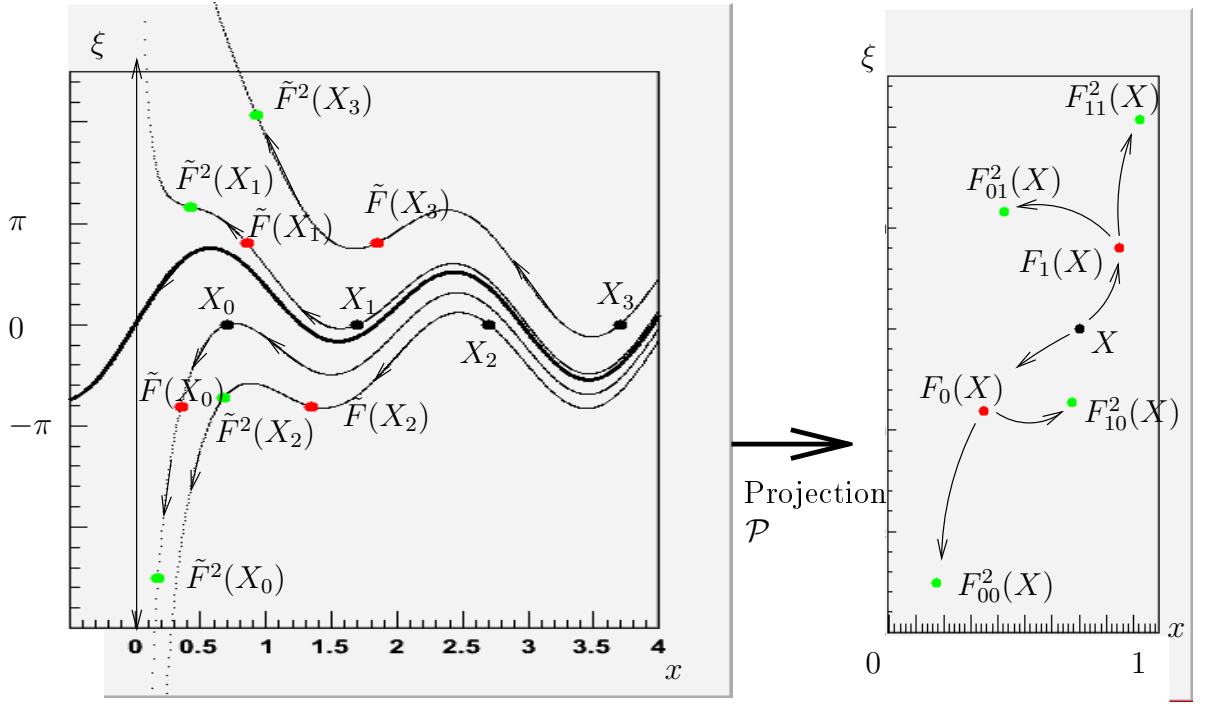
\par\end{centering}

\caption{\label{fig:Relation-dyn}This picture shows how the dynamics of a
point $X=\mathcal{P}\left(X_{0}\right)\in S^{1}\times\mathbb{R}$
under the map $F$ is related by Eq.(\ref{eq:dyna_R2_Cyl}) to the
dynamics of its lifted images $X_{k}=X_{0}+\left(k,0\right)$ under
$\tilde{F}$ on the cover $\mathbb{R}^{2}$.}

\end{figure}

One easily establish the following properties of the map $\tilde{F}$,
illustrated on figure \ref{fig:separatrice}:
\begin{itemize}
\item The point \[
I:=\left(0,\xi_{I}\right):=\left(0,-\frac{\tau'\left(0\right)}{\left(E'\left(0\right)-1\right)}\right)\]
 is the unique \textbf{fixed point} of $\tilde{F}$. It is hyperbolic
with \textbf{unstable manifold} \[
W_{u}=\left\{ \left(0,\xi\right),\xi\in\mathbb{R}\right\} \]
and \textbf{stable manifold}\begin{equation}
W_{s}=\left\{ \left(x,S\left(x\right)\right),x\in\mathbb{R}\right\} \label{eq:def_Ws}\end{equation}
where the $C^{\infty}$ function $S\left(x\right)$ is defined by
the following\textbf{ co-homological equation}, deduced directly from
(\ref{eq:F_lifted})\[
S\left(E^{-1}\left(x\right)\right)=E'\left(E^{-1}\left(x\right)\right)S\left(x\right)+\tau'\left(E^{-1}\left(x\right)\right),\qquad S\left(0\right)=\xi_{I}.\]
The first equation can be written as \[
S\left(x\right)=\frac{1}{E'\left(E^{-1}\left(x\right)\right)}\left(S\left(E^{-1}\left(x\right)\right)-\tau'\left(E^{-1}\left(x\right)\right)\right)\]
 and recursively we deduce that\begin{equation}
S\left(x\right)=-\sum_{p=1}^{\infty}\frac{1}{E'^{\left(-p\right)}\left(x\right)}\tau'\left(E^{\left(-p\right)}x\right)\label{eq:expression_S}\end{equation}
where \begin{equation}
x_{\left(-p\right)}:=E^{\left(-p\right)}\left(x\right):=\left(\underbrace{E^{-1}\circ\ldots\circ E^{-1}}_{p}\right)\left(x\right)\label{eq:def_x-p}\end{equation}
 and \begin{equation}
E'^{\left(-p\right)}\left(x\right):=E'\left(x_{-p}\right)\ldots E'\left(x_{-2}\right)E'\left(x_{-1}\right)\label{eq:Ep}\end{equation}
is the product of derivatives. In the case of $E\left(x\right)=2x$,
one gets simply \begin{equation}
S\left(x\right)=-\sum_{p=1}^{\infty}\frac{1}{2^{p}}\tau'\left(\frac{x}{2^{p}}\right)\label{eq:S_simple}\end{equation}

\item If $\mathcal{P}:\left(x,\xi\right)\in\mathbb{R}^{2}\rightarrow\left(x\, mod\,1,\xi\right)\in S^{1}\times\mathbb{R}\equiv T^{*}S^{1}$
denotes the projection, then \textbf{trapped set} $K$, defined in
Eq.(\ref{eq:def_K}) is obtained by wrapping the stable manifold around
the cylinder and taking the closure:\begin{equation}
K=\overline{\mathcal{P}\left(W_{s}\right)}\label{eq:K_P_Ws}\end{equation}
Compare Figures \ref{fig:separatrice} and Figure \ref{fig:The-trapped-set}.
\item If $X_{0}=\left(x_{0},\xi_{0}\right)\in\mathbb{R}^{2}$ is an initial
point on the plane, and $\mathcal{P}\left(X_{0}\right)$ denotes its
image on the cylinder, then at time $n\in\mathbb{N}$, the $k^{n}$
evolutions of the point $\mathcal{P}\left(X_{0}\right)$ under the
map $F^{n}$ are the images of the evolutions $\tilde{F}^{n}\left(X_{k}\right)$
of the translated points $X_{p}=X_{0}+\left(p,0\right)$, with $p=0\rightarrow k^{n}-1$:\[
F^{n}\left(\mathcal{P}\left(X_{0}\right)\right)=\left\{ \mathcal{P}\left(\tilde{F}^{n}\left(X_{p}\right)\right),\, X_{p}=X_{0}+\left(p,0\right),\quad p\in[0,k^{n}[\right\} \]
and more precisely, using notation of Eq.(\ref{def:trajectory}) for
these points, one has the relation: \begin{equation}
F_{\varepsilon}^{n}\left(\mathcal{P}\left(X_{0}\right)\right)=\mathcal{P}\left(\tilde{F}^{n}\left(X_{p}\right)\right)\label{eq:dyna_R2_Cyl}\end{equation}
where $\varepsilon$ is the number $p$ written in base $k$:\[
\varepsilon=\varepsilon_{n-1}\ldots\varepsilon_{1}\varepsilon_{0}=p_{base\, k}\in\left\{ 0,\ldots,k-1\right\} ^{n}\]
Figure \ref{fig:Relation-dyn} illustrates this correspondence.
\item For an initial point $X_{0}=\left(x_{0},\xi_{0}\right)\in\mathbb{R}^{2}$,
then $X_{n}=\left(x_{n},\xi_{n}\right)=\tilde{F}^{n}\left(X_{0}\right)$
satisfies\begin{equation}
x_{n}=E^{\left(-n\right)}\left(x_{0}\right),\qquad\xi_{n}-S\left(x_{n}\right)=\left(E'^{\left(-n\right)}\left(x_{0}\right)\right)\left(\xi_{0}-S\left(x_{0}\right)\right)\label{eq:trajectoire_R2}\end{equation}
with $E^{\left(-n\right)}\left(x\right)$, $E'^{\left(-n\right)}\left(x\right)$
given by (\ref{eq:def_x-p}), (\ref{eq:Ep}). Hence\[
\left|\xi_{n}-S\left(x_{n}\right)\right|\geq E_{min}^{n}\left|\xi_{0}-S\left(x_{0}\right)\right|\]
This last inequality describes how fast the trajectories above or
below the separatrix $W_{s}$ escape towards infinity on Figure \ref{fig:separatrice}.
\end{itemize}

\subsection{Partially captive property}

Here we rephrase the property of partial captivity, Definition \ref{def:partially_captive},
in terms of a property on the separatrix function $S\left(x\right)$
defined in Eq.(\ref{eq:def_Ws}) and given in Eq. (\ref{eq:expression_S}).

For simplicity, we consider from now on the simple model with a linear
expanding map $E\left(x\right)=kx$.

\vspace{0.cm}\begin{center}{\color{blue}\fbox{\color{black}\parbox{15cm}{
\begin{prop}

For $n\in\mathbb{N}$, and $\tilde{R}>0$, let

\begin{equation}
\tilde{\mathcal{N}}_{\tilde{R}}\left(n\right)=\max_{\left(x,\xi\right)\in\mathbb{R}^{2}}\sharp\left\{ p\in[0,k^{n}[,\quad\left|\xi-S\left(x+p\right)\right|\leq\frac{\tilde{R}}{k^{n}}\right\} \label{eq:prop_S}\end{equation}

Then the map $F$ is \textbf{partially captive} (see definition page
\pageref{def:partially_captive}) if and only if \begin{equation}
\lim_{n\rightarrow\infty}\frac{\log\tilde{\mathcal{N}}_{\tilde{R}}\left(n\right)}{n}=0\label{eq:partially_captive_bis}\end{equation}
for $\tilde{R}$ large enough.

\end{prop}
}}}\end{center}\vspace{0.cm}
\begin{proof}
From (\ref{eq:def_Nn}) and using (\ref{eq:dyna_R2_Cyl}),(\ref{eq:trajectoire_R2})
one gets\begin{eqnarray*}
\mathcal{N}\left(n\right) & = & \max_{\left(x,\xi\right)}\sharp\left\{ F_{\varepsilon}^{n}\left(\left(x,\xi\right)\right)\in\mathcal{Z},\qquad\varepsilon\in\left\{ 0,\ldots,k-1\right\} ^{n}\right\} \\
 & = & \max_{\left(x,\xi\right)}\sharp\left\{ \left|\xi_{n,p}\right|\leq R,\qquad p\in[0,k^{n}[\right\} \end{eqnarray*}
with $\left(x_{n,p},\xi_{n,p}\right):=\tilde{F}^{n}\left(\left(x+p,\xi\right)\right)$
given by $x_{n,p}=\frac{x+p}{k^{n}}$ and $\xi_{n,p}=S\left(x_{n,p}\right)+k^{n}\left(\xi-S\left(x+p\right)\right)$.
Therefore\[
\left|\xi_{n,p}\right|\leq R\Leftrightarrow\left|\frac{S\left(x_{n,p}\right)}{k^{n}}+\xi-S\left(x+p\right)\right|\leq\frac{R}{k^{n}}\]
 But $S$ is a bound function on $\mathbb{R}$, $\left|S\left(x\right)\right|\leq S_{max}$.
Therefore $\left|\frac{S\left(x_{n,p}\right)}{k^{n}}+\xi-S\left(x+p\right)\right|\leq\frac{R}{k^{n}}\Longrightarrow\left|\xi-S\left(x+p\right)\right|\leq\frac{\tilde{R}}{k^{n}}$
with some $\tilde{R}>0$, and conversely. This implies that (\ref{eq:partially_captive_bis})
is equivalent to (\ref{eq:def_partially_captive}).
\end{proof}

\subsection{Fractal aspect of the Trapped set}

The characterization Eq.(\ref{eq:prop_S}) concerns the discrete set
of points $S\left(x+m\right),\,\, m\in\mathbb{Z}$. From Eq.(\ref{eq:K_P_Ws})
these points are the slice of the trapped set $K=\cup_{x\in S^{1}}K_{x}$:\[
K_{x}=\overline{\left\{ S\left(x+m\right),\, m\in\mathbb{Z}\right\} }\]

For simplicity, we consider from now on the simple model with $E\left(x\right)=2x$,
and $\tau\left(x\right)=\cos\left(2\pi x\right)$.

\begin{figure}[h]
\begin{centering}
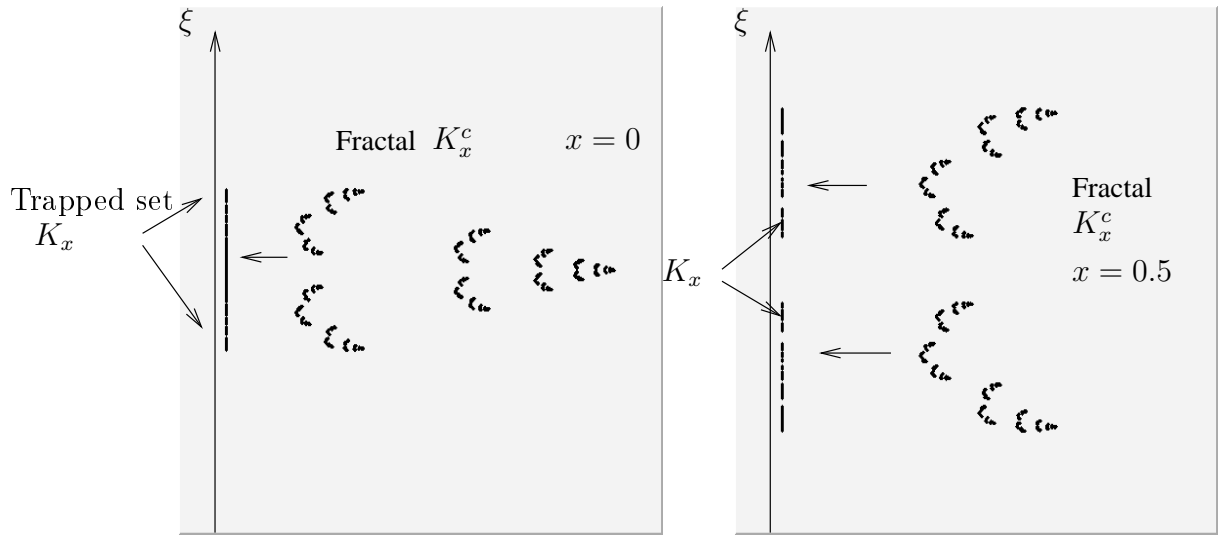
\par\end{centering}

\caption{\label{fig:fractal} This picture represents the set $K_{x}^{c}\subset\mathbb{C}$
defined by Eq.(\ref{eq:def_Kxc}). The trapped set $K$ at position
$x$ is obtained by the projection on the imaginary axis $K_{x}=\Im\left(K_{x}^{c}\right)$.
On the web page \protect\url{http://www-fourier.ujf-grenoble.fr/~faure/articles}
one can observe the motion of the fractal $K_{x}^{c}$ as $x\in\mathbb{R}$
increases smoothly.}

\end{figure}

From Eq.(\ref{eq:S_simple}), these points are given by\[
S\left(x+m\right)=\sum_{p=1}^{\infty}\frac{2\pi}{2^{p}}\sin\left(\frac{2\pi}{2^{p}}\left(x+m\right)\right)=\Im\left(\sum_{p=1}^{\infty}\frac{2\pi}{2^{p}}\exp\left(\frac{i2\pi}{2^{p}}\left(x+m\right)\right)\right)\]
Therefore the slice $K_{x}$ is the projection on the imaginary axis
of the following set:\begin{equation}
K_{x}^{c}=\overline{\left\{ S^{c}\left(x+m\right),m\in\mathbb{Z}\right\} }\label{eq:def_Kxc}\end{equation}
\[
S^{c}\left(x+m\right)=\sum_{p=1}^{\infty}\frac{2\pi}{2^{p}}e^{i2\pi\frac{1}{2^{p}}\left(x+m\right)}\]
On Figure \ref{fig:fractal} we observe that $K_{x}^{c}$ is \textbf{a
fractal set.} Compare Figure \ref{fig:fractal} with Figure \ref{fig:The-trapped-set}.

\bibliographystyle{alpha}
\bibliography{/home/faure/articles/articles}

\end{document}